\documentclass[10pt]{article}

\usepackage{hyperref}

\hypersetup{
            pdffitwindow=false,
            plainpages=false,
            pdfpagemode=UseOutlines,
            pdfpagelayout=SinglePage,
            colorlinks=true,
            linkcolor=blue,
            urlcolor=blue!30!black,
            citecolor=green!50!black,
}

\usepackage[ngerman,english]{babel} 
\usepackage{a4wide}
\usepackage{amsmath}
\usepackage{amssymb}
\usepackage{amsthm}
\usepackage[title]{appendix}
\usepackage{bm}
\usepackage{graphicx}
\usepackage{caption}
\usepackage{subcaption}
\usepackage[latin1]{inputenc}
\usepackage{enumerate}
\usepackage{ifthen}
\usepackage{verbatim}
\usepackage{color}
\usepackage{float}
\usepackage{pgfplots}
\usepackage{multicol}
\usepackage{mathtools}
\usepackage{float}
\usepackage{bm}
\usepackage{bbm}
\usepackage{dsfont}
\usepackage{placeins}

\newcommand{\R}{\mathbb{R}}

\newcommand{\ep}{\varepsilon}

\newcommand{\op}[1]{\mathcal{#1}} 
\newcommand{\set}[1]{\mathbb{#1}}

\DeclareMathOperator{\supp}{supp}

\definecolor{boxback}{gray}{0.95}

\theoremstyle{plain}
\newtheorem{theorem}{Theorem}[section]
\newtheorem{lemma}[theorem]{Lemma}
\newtheorem*{remark*}{Remark}
\newtheorem{remark}[theorem]{Remark}

\theoremstyle{definition}

\usepgfplotslibrary{external}
\usepgfplotslibrary{groupplots}
\usepgfplotslibrary{fillbetween}
\usetikzlibrary{pgfplots.groupplots}
\usetikzlibrary{matrix}
\usetikzlibrary{plotmarks}

\begin{document}

\title{Diffusion maps tailored to arbitrary non-degenerate It\^o processes}


\author{Ralf Banisch \and Zofia Trstanova \and Andreas Bittracher \and Stefan Klus \and P\'eter Koltai}



\maketitle


\begin{abstract}
We present two generalizations of the popular diffusion maps algorithm. The first generalization replaces the drift term in diffusion maps, which is the gradient of the sampling density, with the gradient of an arbitrary density of interest which is known up to a normalization constant. The second generalization allows for a diffusion map type approximation of the forward and backward generators of general It\^o diffusions with given drift and diffusion coefficients. We use the local kernels introduced by Berry and Sauer, but allow for arbitrary sampling densities. We provide numerical illustrations to demonstrate that this opens up many new applications for diffusion maps as a tool to organize point cloud data, including biased or corrupted samples, dimension reduction for dynamical systems, detection of almost invariant regions in flow fields, and importance sampling.
\end{abstract}

\section{Introduction}

The analysis of large data sets in Euclidean space is a topic of active research, with methods such as dimension reduction and manifold learning trying to identify intrinsic characteristics of the data. The central goals are, in essence, to compress the data from an originally large number of variables into a low-dimensional description involving only a few variables, and to find meaningful structure in the data in order to gain insight and understanding of the process that generated it.

A common starting point is principal component analysis (PCA), which tries to understand data by performing a singular value decomposition of the data matrix and then taking only the singular vectors corresponding to the largest singular values as coordinates. This will be successful if the data lies on a linear subspace, but nonlinear structures will not be uncovered. Kernel PCA \cite{Ham2004, Schoelkopf1997} tries to understand nonlinear features in the data by mapping it to a usually high-dimensional feature space first, this is done implicitly by constructing a symmetric positive definite matrix $K$ of inner products. An alternative approach assumes that the data points lie on a manifold $\op{M}$ and attempts to encode structure through differential operators on $\op{M}$. The first step is typically the construction of a neighborhood graph with similarity weights derived from a kernel function. Then a graph Laplacian is constructed and the dominant eigenvectors are used to organize the data \cite{von2007tutorial}.

Belkin and Niyogi \cite{Belkin2003} and Coifman and Lafon \cite{CoLa06} use radially symmetric kernels and construct a graph Laplacian that estimates the Laplace operator on the manifold $\op{M}$. The connection between graph Laplacians based on radially symmetric kernels and the Laplace operator on the sampled manifold was further understood in \cite{hein2005graphs}. Berry and Sauer \cite{BeSa16} extended this approach to a much larger class of anisotropic local kernels and showed a connection between the moment expansion of the kernel and the differential operators that can be approximated. This allows the approximation of the forward and backward Fokker--Planck operators associated with a large class of It\^o diffusions, provided a kernel with prescribed first and second moment is used.

Often the data will not sample the manifold uniformly, and care must be taken to distinguish effects that are due to local variations in the sampling density $q$ from intrinsic characteristics of the manifold. The approach in \cite{Belkin2003} and, as we show in this article, \cite{BeSa16} produces unwanted drift terms if the sampling is not uniform which bias the approximation of the considered differential operators. If, as is often the case, the sampling density is unknown, then this bias cannot be controlled. A key property of the diffusion maps construction \cite{CoLa06} is that by a suitable modification of the kernel the bias induced by local variations in $q$ can be completely removed, which leads to asymptotically unbiased estimators of the Laplacian and the heat kernel describing diffusion on~$\op{M}$.

A second key property of diffusion maps is that with a slightly different modification of the kernel, the diffusive terms encoded in the kernel and the bias from the sampling density $q$ can be combined in order to approximate the generator of a Markov process in the case of a gradient flow \cite{Nadler2006, Nadler2008, singer2009detecting, szlam2008regularization} where the gradient terms are determined by the gradient of $q$. This interesting property allows the approximation of the generator of a dynamical process based purely on the geometry of the data. The differential operator approximated in this case depends explicitly on $q$. The data does not have to be generated by the dynamical process, but it has to be sampled from $q$ in a controlled fashion, which can be challenging for many densities of interest. Nevertheless this view has led to many applications where Markov processes with a gradient flow are interesting \cite{Nadler2006, Nadler2008, singer2009detecting, rohrdanz2011determination, Noe2015, Noe2016}.

In this article, we introduce two extensions of the diffusion maps construction. Both extensions are asymptotically unbiased in the sense that the limiting differential operators do not depend on $q$. The first extension, target measure diffusion maps, generalizes the construction in \cite{CoLa06, Nadler2006, Nadler2008, singer2009detecting, szlam2008regularization} and allows to approximate the generator of a gradient flow Markov process where the gradient terms are not slaved to the sampling density, but are derived from a target probability measure that is chosen by the user and known up to a normalization constant. This allows us to extend diffusion maps to cases where sampling the measure of interest is either impossible or too difficult, and one has to work with biased samples instead. Our approach is related to \cite{Berry2016}, which uses variable bandwidth kernels and requires knowledge of the intrinsic dimension of $\op{M}$.

The second extension, local kernel diffusion maps, allows to approximate the forward and backward Fokker--Planck operators of a large class of It\^o diffusions on $\op{M}$. We use the anisotropic local kernels introduced in \cite{BeSa16}, but we also estimate the sampling density $q$ by using a second, radially symmetric kernel and are able to undo its influence. This generalizes \cite{Nadler2006, Nadler2008, singer2009detecting, szlam2008regularization} to non-gradient flows and anisotropic noise, and the results from \cite{BeSa16} to non-uniform sampling densities. We also discuss how the drift and diffusion coefficients can be estimated from data if they are not known analytically. This extends diffusion maps to approximate the generators of a much larger class of dynamical processes.

In Section \ref{sec:diffmap} we summarize diffusion maps for radial kernels as found in \cite{CoLa06, Nadler2006, Nadler2008, singer2009detecting, SINGER2006} and discuss the results and developments relevant for this article. In Section \ref{sec:anisotropic} we discuss the anisotropic local kernels introduced in \cite{BeSa16} and show that they lead to biased estimators in the case of nonuniform sampling. Our main results can be found in Sections \ref{sec:TMDmap} and \ref{sec:localkernel_unbiased}, where we introduce target measure and local kernel diffusion maps and state the respective convergence results. Many numerical examples that showcase applications are discussed in Section \ref{sec:numerics} and proofs can be found in Section \ref{sec:proofs}.

\section{Diffusion maps for isotropic kernels}\label{sec:diffmap}

In this article, we assume that we are given data in the form of points $\set{D}^{(m)} := \{x_1,x_2,\ldots,x_m\}\subset\R^N$ with $N>0$, which lie on a compact $d$-dimensional differentiable submanifold $\op{M}\subset \R^N$ and are sampled according to a density $q(x)$. We assume that both $q$ and $\op{M}$ are unknown. The dimension $d$ of $\op{M}$ is also assumed to be unknown, but we will have the case $d\ll N$ in mind. The data set $\set{D}^{(m)}$ is used to represent functions $f$ as vectors $[f] = (f(x_1), \ldots, f(x_m))^T$ and operators $\op{A}$ as $m\times m$ matrices $A$. The $i$th component of a data vector is denoted by $[f]_i = f(x_i)$, and $(A[f])_i$ is the $i$th component of the matrix-vector product $A[f]$. We denote the transpose of a matrix $A$ by~$A^T$. The constant function $f(x) \equiv 1$ is denoted by $\mathbf{1}$, the corresponding constant data vector by~$[\mathbf{1}]$.

The main idea underlying diffusion maps~\cite{CoLa06} is to uncover the geometric structure of $\op{M}$ from the data~$\set{D}^{(m)}$ by constructing an $m\times m$ matrix that approximates a differential operator. This differential operator, and in particular its dominant eigenfunctions, encodes all relevant geometric features of $\op{M}$. 

We briefly summarize the algorithm and most relevant results in \cite{CoLa06, SINGER2006}. The construction starts with an isotropic kernel $k_\ep(x, y) = h( \|x - y\|^2 / \ep)$ where $h$ is an exponentially decaying function, $\ep > 0$ is a scale parameter and $\|\cdot\|$ the Euclidean norm in $\R^N$. The most common choice is
\begin{equation}
k_{\ep}(x,y) = \exp\left(-(4\ep)^{-1}\|x-y\|^2\right).
\label{eq:kernel}
\end{equation}
Evaluation of $k_\ep$ on $\set{D}^{(m)}$ leads to a $m\times m$ kernel matrix $K_\ep$, which is normalized several times to give a matrix $L_{\ep,\alpha}$ that can be interpreted as the generator of a Markov chain on the data itself. Here, $\alpha \in \R$ is a parameter that can be chosen during the construction. In the limit $m\rightarrow \infty$ and $\ep \rightarrow 0$, the matrix $L_{\ep, \alpha}$ approximates the Kolmogorov operator
\begin{equation}
\label{eq:gen_diffmaps_vanilla}
\op{L}f = \Delta f + (2-2\alpha)\nabla f\cdot \frac{\nabla q}{q},
\end{equation}
where $\Delta$ is the Laplace--Beltrami operator on $\op{M}$ and $\nabla$ the gradient operator. The approximation is in the sense that for sufficiently smooth functions $f$ and any point $x_k \in \set{D}^{(m)}$, $(L_{\ep,\alpha} [f])_k \rightarrow \op{L}f(x_k)$ as~$m\to\infty$. Moreover, if $q$ is uniform then $L_{\ep,0}[f]_k = \op{L} f(x_k) + \mathcal{O}(\ep,m^{-1/2} \ep^{-1/2-d/4})$ with high probability \cite{SINGER2006}. Consequently, eigenfunctions of $L_{\ep,\alpha}$ approximate eigenfunctions of $\op{L}$. If $\op{M}$ is bounded, then eigenfunctions of $L_{\ep, \alpha}$ approximate solutions to the eigenproblem of $\op{L}$ with von Neumann boundary conditions on $\partial \op{M}$. Spectral convergence of $L_{\ep,0}$ was shown in \cite{vonluxburg2008}.

The eigenvalues of $L_{\ep,\alpha}$ are real and nonpositive and can be ordered as $0 = \lambda_0 > \lambda_1 \geq \lambda_2 \geq \ldots \geq \lambda_m$. The largest eigenvalue is always $\lambda_0 = 0$, and under mild assumptions it has multiplicity one. Of particular interest are the \emph{dominant} eigenvalues, i.e.~the smallest in magnitude but nonzero eigenvalues. These eigenvalues give the slowest time scales of the diffuse dynamics generated by \eqref{eq:gen_diffmaps_vanilla}, and the corresponding eigenfunctions allow for a structure preserving embedding $\Psi$ of $\set{D}^{(m)}$ into a lower-dimensional space. This \emph{diffusion map} is a way of representing the geometry of the data.

Perhaps the two most well-known special cases of this construction are (a) $\alpha = 1$, in which case $\op{L} = \Delta$ is the Laplace--Beltrami operator and dominant eigenfunctions are coordinates on $\op{M}$ and (b) $\alpha = 1/2$ and $q \propto \exp(-\beta U)$ is the Boltzmann density with parameter $\beta > 0$, in which case
\begin{equation}
\label{eq:diffmap_standard}
\op{L}f = \Delta f - \beta\nabla U\cdot \nabla f,
\end{equation}
which is $\beta$ times the backward Kolmogorov operator associated with Brownian motion at temperature $\beta^{-1}$ in the potential energy landscape given by the energy function $U: \op{M}\rightarrow \R$. This powerful property of diffusion maps allows the computation of the eigenfunctions of such a dynamics based purely on geometric data produced by sampling from the Boltzmann distribution. It is well known that the dominant eigenvalues of the operator \eqref{eq:diffmap_standard} correspond to the equilibration rates of the slow processes \cite{FrJuKo13, Web11, BJK15, SchSa13},  and the associated dominant eigenfunctions allow for dimensionality reduction while retaining the slow dominant time scales \cite{Bittracher2017}. Plateaus in the range of the dominant eigenfunctions indicate metastable subsets of phase space and can be identified by a simple cluster analysis \cite{DeJu99, SchueFi99} or more elaborate means~\cite{DeWe03}. These properties have led to several applications in the Molecular Dynamics context \cite{Nadler2008, Nadler2006, rohrdanz2011determination, Noe2015, Noe2016}. We emphasize that it does not matter how the samples were generated as long as they are distributed according to  $q \propto \exp(-\beta U)$.

To appreciate the limitations of diffusion maps, consider an It\^o stochastic differential equation in $\R^N$ of the form
\begin{equation}
dX_t = b(X_t)\,dt + \sigma(X_t)\,dW_t
\label{eq:SDE}
\end{equation}
with drift $b: \R^N \rightarrow \R^N$, diffusion coefficient $\sigma: \R^N \rightarrow \R^{N\times M}$ and $W_t$ being $M$-dimensional Brownian motion. It is well known that the generator $\op{L}$ of the dynamics \eqref{eq:SDE} is given by the backward Kolmogorov operator
\begin{equation}
\label{eq:generator_general}
\op{L}f = b\cdot \nabla f + A_{ij} \nabla_i \nabla_j f
\end{equation}
with diffusion matrix $A = \frac{1}{2} \sigma\sigma^T$. Here, $\nabla_i$ is the $i$th component of the operator $\nabla$ and we use the Einstein summation convention. Dot products are taken with respect to the metric inherited from the ambient space $\R^N$. The family of operators \eqref{eq:gen_diffmaps_vanilla} that can be approximated by diffusion maps corresponds to the special case $A = I$ and $b = (2-2\alpha) q^{-1}\nabla q$. Therefore, if we wish to use diffusion maps to approximate generators of It\^o processes \eqref{eq:SDE}, there are two key limitations:
\begin{itemize}
\item Isotropic noise imposed by $A = I$. This is a direct consequence of the use of an isotropic kernel.
\item The drift $b$ must either be zero (which is the case $\op{L} = \Delta$), or it must be related to the sampling density $q$ via $b \propto \nabla\log q$. This means that if we wish to ensure that $b = -\beta \nabla U$ with a prescribed energy function $U : \R^N \rightarrow \R$, then $q$ must be the Boltzmann density $q = \exp(-\beta U)$.
\end{itemize}

In what follows, we will lift these limitations: We will allow for arbitrary sampling densities, and we will construct a variant of diffusion maps that can reconstruct the generator of an arbitrary SDE of the form \eqref{eq:SDE}.

\section{Target measure diffusion maps}\label{sec:TMDmap}

Let $d\pi$ be a probability measure on $\R^N$ whose density $\pi$ is known only up to a normalization constant, which is the case in many applications such as Markov Chain Monte Carlo. We do not assume that $\pi$ is related to the unknown sampling density $q$. Suppose that we wish to approximate the differential operator
\begin{equation}
\op{L} f = \Delta f + \nabla (\log \pi) \cdot \nabla f,
\label{eq:backward_gen_pi}
\end{equation}
which is the generator of the It\^o diffusion
\begin{equation}
\label{eq:SDE_ergodic}
dX_t = \nabla(\log(\pi(X_t)) + \sqrt{2}dW_t.
\end{equation}
Under mild assumptions on $\pi$, the dynamics \eqref{eq:SDE_ergodic} is ergodic with invariant measure $d\pi$~\cite{meyn2012markov, roberts1996exponential}. The adjoint of $\op{L}$ is given by
\[
\op{L}^* = \Delta f - \nabla\cdot \left( f \, \nabla(\log \pi)\right)
\]
and satisfies $\op{L}^* \pi = 0$. Often $\pi$ is a Boltzmann density of the form $\pi = Z^{-1}\exp(-\beta U)$, where $\beta>0$ is an inverse temperature parameter, $U : \R^N \rightarrow \R$ an energy function and $Z$ the unknown normalization constant. In this case, $\op{L}f = \Delta f - \beta \nabla U \cdot \nabla f$ is the operator in \eqref{eq:diffmap_standard}.

In general, obtaining samples directly from the target distribution $\pi$ is difficult for several reasons: First, the device that generates the data $\set{D}^{(m)}$ could have a bias, so that it samples not $\pi$, but some distribution close to $\pi$. This is already the case if the dynamics \eqref{eq:SDE_ergodic} is discretized in time. One may remove the bias by Metropolization~\cite{hastings1970monte,roberts1996exponential}, but this renders sampling more expensive since proposal steps might be rejected. Second, we may wish to introduce a bias intentionally, in the sense of importance sampling. If, for example, $\beta^{-1}$ is small compared to the energy barriers in $U$, then the dynamics \eqref{eq:SDE_ergodic} is known to suffer from metastability issues and will therefore converge exponentially slowly to its invariant measure $d\pi$~\cite{stoltz2010free}. In order to accelerate convergence, we may wish to sample at a higher temperature so that $q \propto \exp(-\tilde \beta U)$ with $\tilde \beta < \beta$. Third, it is often only possible to generate correlated samples, i.e.~by discretizing \eqref{eq:SDE_ergodic} in time. These correlations can have long decay times, and it is possible that all our samples $x_i$ have strong correlations with the initial sample $x_1$.

Motivated by this, we generalize diffusion maps so that the requirement $q=\pi$ is no longer necessary: We introduce the \emph{target measure diffusion map} (TMDmap) algorithm which constructs a matrix $L_{\ep, \pi}$ on the data that approximates the differential operator \eqref{eq:backward_gen_pi}. We use the isotropic kernel \eqref{eq:kernel} and assume that the target density $\pi$ can be evaluated on the data up to a normalization constant.
\begin{center}
\colorbox{boxback}{
  \begin{minipage}[t]{\textwidth}
\begin{center}
  \begin{minipage}[t]{0.95\textwidth}
\vspace{0.5cm}
\center\textbf{Construction of TMDmap}
\begin{enumerate}[(i)]
\item Construct the $m\times m$ kernel matrix $K_\ep$ with components $(K_\ep)_{ij} = k_\ep(x_i, x_j)$ and the kernel density estimate
\[
q_\ep(x_i) = \sum_{j=1}^m (K_\ep)_{ij}.
\]
\item Form the diagonal matrix $D_{\ep,\pi}$ with components $(D_{\ep,\pi})_{ii} = \pi^{1/2}(x_i)\,q^{-1}_\ep(x_i)$ and right-normalize the kernel matrix with $D_{\ep,\pi}$:
\[
K_{\ep,\pi} = K_\ep D_{\ep,\pi}.
\]
\item Let $\tilde D_{\ep,\pi}$ be the diagonal matrix of row sums of $K_{\ep,\pi}$, that is,
\[
(\tilde D_{\ep,\pi})_{ii} = (K_{\ep,\pi} [\mathbf{1}])_i = \sum_{j=1}^m (K_{\ep, \pi})_{ij}.
\]
\item Build the TMDmap matrix
\begin{equation}
\label{eq:TMDmap}
L_{\ep, \pi} = \ep^{-1}\left( \tilde D_{\ep, \pi}^{-1} K_{\ep,\pi} - I\right).
\end{equation}
\end{enumerate}
  \end{minipage}
\vspace{0.5cm}
\end{center}
 \end{minipage}
}
\end{center}
The idea of TMDmap is to normalize the kernel twice: The right normalization of $K_\ep$ with the diagonal matrix $D_{\ep,\pi}$ in step (ii) cancels unwanted drift terms produced by the sampling density $q$, which is estimated via the kernel density estimator $q_\ep$. It also introduces the desired drift terms coming from $\pi$. The left normalization of $K_{\ep,\pi}$ by $\tilde D_{\ep,\pi}$ in step (iii) ensures that $L_{\ep,\pi}$ conserves probability and can thus be interpreted as the generator of a Markov process on the data.

The following theorem establishes the relation between the matrix $L_{\ep,\pi}$ and the differential operator \eqref{eq:backward_gen_pi}. The pointwise convergence established here is completely equivalent to the analogous diffusion map results established in \cite{CoLa06}. The proof can be found in Appendix \ref{sec:proof_TMDmap}.

\begin{theorem}[Convergence of TMDmap]
\label{thm:target_diffmap}
Let $L_{\ep,\pi}$ be defined via \eqref{eq:TMDmap} and let $f : \op{M} \rightarrow \R$ be smooth. Then in the limit $m\rightarrow \infty$ and $\ep \rightarrow 0$ and for every $x_k \in \set{D}^{(m)}$,
\begin{equation*}
(L_{\ep,\pi}[f])_k \rightarrow \Delta f(x_k) + \nabla (\log \pi) \cdot \nabla f(x_k).
\end{equation*}
\end{theorem}
Theorem \ref{thm:target_diffmap} establishes pointwise convergence of $L_{\ep,\pi}[f]$ to the action of the differential operator $\op{L}$ in \eqref{eq:backward_gen_pi} on the function $f$, evaluated at the data points. Typically, eigenfunctions of $\op{L}$ are the objects of interest. As in standard diffusion maps \cite{CoLa06}, solutions of the matrix eigenvalue problem $L_{\ep, \pi} [\psi] = \lambda [\psi]$ can be used to approximate solutions of the eigenvalue problem
\begin{equation}
\label{eq:TMDmap_eval}
\op{L}\psi (x) = \lambda\psi(x) \quad \forall x\in \op{M} = \supp(q).
\end{equation}
If $\op{M}$ is bounded, then von Neumann boundary conditions are imposed in \eqref{eq:TMDmap_eval}. 

We first note that TMDmap is in fact a generalization of standard diffusion maps: In the special case where the sampling density $q$ is equal to $\pi$, the result in Theorem \ref{thm:target_diffmap} is consistent with standard diffusion maps with parameter $\alpha = 1/2$. Unlike standard diffusion maps however, the limiting operator $\op{L}$ in Theorem \ref{thm:target_diffmap} is independent of $q$. The limiting eigenvalue problem of $L_{\ep,\pi}$ is also independent of $q$, but it does depend on $\op{M} = \supp(q)$ in that it fixes the domain of \eqref{eq:TMDmap_eval}. This is to be expected; we cannot hope to approximate $\op{L}$ in regions we will never draw samples from. TMDmap is probably most useful in situations where $\supp(q) \supseteq \supp(\pi)$, which includes the biased sampling and sampling at higher temperatures examples discussed above. If $\supp(q) \subset \supp(\pi)$, then eigenvectors of $L_{\ep,\pi}$ will approximate eigenfunctions of $\op{L}$ in the restricted domain $\op{M}$.

We end the section by comparing TMDmap with an alternative approach to approximate the operator \eqref{eq:backward_gen_pi} with $\pi \propto \exp(-\beta U)$ proposed in \cite{Berry2016}. The idea is to use the variable bandwidth kernel
\[
K^S_\ep(x,y) = h\left(\frac{\|x-y\|^2}{\ep \rho(x)\rho(y)}\right)
\]
where $h$ is an exponentially decaying function and $\rho$ is a bandwidth function. A diffusion map type construction detailed in \cite{Berry2016} with the choice $\rho(x) = \exp(-\beta U / (d+2))$, where $d = \mbox{dim}(\op{M})$, will also approximate \eqref{eq:backward_gen_pi}. A limitation is that the intrinsic dimension of $\op{M}$ has to be explicitly known in this case, which is not necessary for TMDmap. On the other hand, variable bandwidth kernels may show less sensitivity with respect to $\ep$ and may have more favourable convergence properties, especially for unbounded $\op{M}$.

\begin{remark}
It is worth noting that
\[
D_{\ep,\pi}\tilde D_{\ep,\pi}L_{\ep,\pi} = \ep^{-1}(D_{\ep,\pi}K_{\ep,\pi}D_{\ep,\pi} - D_{\ep,\pi}\tilde D_{\ep,\pi})
\]
is a symmetric matrix. Therefore, left eigenvectors $\phi_n^T L_{\ep,\pi} = \lambda_n \phi_n^T$ and right eigenvectors $L_{\ep,\pi}\psi_n = \lambda_n\psi_n$ are related via $\phi_n = D_{\ep,\pi}\tilde D_{\ep,\pi}\psi_n$. Moreover, the right eigenvector $\psi_0$ corresponding to $\lambda_0 = 0$ is the constant function since $L_{\ep,\pi}$ has row sum zero. The left zero eigenvector $\phi_0$ is thus given by $\phi_0 = \text{diag}(D_{\ep,\pi}\tilde D_{\ep,\pi})$, or
\[
\phi_0 = \frac{\pi^{1/2}(x_i)}{q_\ep(x_i)}\sum_{j=1}^m k_\ep(x_i, x_j) \frac{\pi^{1/2}(x_j)}{q_\ep(x_j)}.
\]
From here, it is easy to see that
\begin{equation}
\label{eq:pi_hat}
\phi_0 = \frac{\pi(x_i)}{q_\ep(x_i)} + \mathcal{O}(\ep),
\end{equation}
so the target density $\pi$ is recovered as the left zero eigenvector $\phi_0$ of $L_{\ep,\pi}$ times the kernel density estimate $q_\ep$. The vector $\phi_0$ contains the weights that are necessary in order to compute $\pi$-averages with our $q$-distributed samples.
\end{remark}

\section{Anisotropic kernels}\label{sec:anisotropic}

Suppose one wishes to extend the family of operators \eqref{eq:diffmap_standard} that can be approximated with diffusion maps to the larger family of forward and backward Kolmogorov operators corresponding to the It\^o process \eqref{eq:SDE} with general---and possibly position dependent---drift $b$ and diffusion matrix~$A$:
\begin{subequations}
\begin{align}
\op{L}f &\; := b \cdot \nabla f + A_{ij}\nabla_i\nabla_j f , \label{eq:backward_gen}\\
 \op{L}^*f &\; = -\nabla\cdot (b f) + \nabla_i \nabla_j (A_{ij} f).\label{eq:forward_gen}
\end{align}
\end{subequations}
The restriction to isotropic diffusion $A=I$ in standard diffusion maps is a consequence of the use of the isotropic kernel \eqref{eq:kernel}. This suggests to modify the kernel. Various attempts have been made to generalize diffusion maps to anisotropic kernels \cite{rohrdanz2011determination, singer2008non, singer2009detecting, dsilva2013nonlinear}. Berry and Sauer~\cite{BeSa16} give a definitive answer in the case that the sampling density $q$ is uniform. In this section, we summarize the results in \cite{BeSa16} for our setting and explicitly include the effect of the sampling density. We then show that unwanted drift terms appear if $q$ is not uniform.

Let $A(x)$ be a matrix-valued function on $\op{M}$ such that each $A(x)$ is a symmetric positive definite $N\times N$ matrix, and let $b(x)$ be a vector-valued function. Berry and Sauer introduce the family of \emph{local kernels}
\begin{equation}
k^{A,b}_{\ep}(x,y) = \exp\left(-\left(4\ep\right)^{-1} \left(x-y+\ep b(x)\right)^{T}A(x)^{-1}\left(x-y+\ep b(x)\right)\right)\,,
\label{eq:localkernel}
\end{equation}
and show that one may approximate $\op{L}$ and $\op{L^*}$ with matrices $L_\ep$ and $L_\ep^*$ defined on the data in the diffusion map sense, provided $q$ is uniform. The construction of the matrices $L_\ep$ and $L_\ep^*$ is as follows. Let $K_\ep^{A,b}$ be the $m\times m$ kernel matrix with components $(K_\ep^{A,b})_{ij} = k^{A,b}_\ep (x_i, x_j)$. Let $D_\ep$ be the diagonal matrix of row sums of $K^{A,b}_\ep$, i.e.
\[
(D_\ep)_{ii} = \sum_{j=1}^m k^{A,b}_\ep(x_i, x_j).
\]
Then the matrices $L_\ep$ and $L^*_\ep$ are constructed as
\begin{subequations}
\begin{align}
L_\ep &\; = \ep^{-1}\left(D_\ep^{-1} K_\ep^{A,b} - I\right)  \label{eq:localkernels_L}\\
L^*_\ep &\; = \ep^{-1}\left((K_\ep^{A,b})^TD_\ep^{-1} - I\right). \label{eq:localkernels_Lstar}
\end{align}
\end{subequations}
One confirms easily that $L^*_\ep = L_\ep^T$ and that $L_\ep$ has row sum zero and positive off diagonal elements, which allows to interpret $L_\ep$ and $L_\ep^*$ as the backward and forward generators of a Markov chain defined on the data. The following theorem, whose proof can be found in Appendix \ref{sec:proof_localkernel}, is a generalization of a result by Berry and Sauer \cite{BeSa16}.

\begin{theorem}
\label{thm:localkernels_q}
Let $A(x)$ be a matrix-valued function on $\op{M}\subset \R^N$ such that each $A(x)$ is a symmetric positive definite $N\times N$ matrix, and let $b(x)$ be a vector-valued function. Let $L_\ep$ and $L_\ep^*$ be defined via \eqref{eq:localkernels_L} and \eqref{eq:localkernels_Lstar} respectively, and let $f : \op{M} \rightarrow \R$ be smooth. Then in the limit $m\rightarrow \infty$ and $\ep \rightarrow 0$ and for every $x_k \in \set{D}^{(m)}$, we obtain
\begin{align*}
(L_\ep [f])_k &\; \rightarrow \op{L}f(x_k) + 2\left[q^{-1}A_{ij} (\nabla_i f) (\nabla_j q)\right] (x_k), \\
(L_\ep^* [f])_k &\; \rightarrow \op{L}^*f(x_k) - \left[q^{-1} f(b\cdot \nabla q +  A_{ij}\nabla_i \nabla_j q)\right] (x_k),
\end{align*}
where $\op{L}$ and $\op{L}^*$ are defined in \eqref{eq:backward_gen} and \eqref{eq:forward_gen}, respectively.
\end{theorem}

If $q$ is uniform, then the unwanted drift terms in Theorem \ref{thm:localkernels_q} are equal to zero and we obtain $(L_\ep [f])_k \rightarrow \op{L}f(x_k)$ and $(L_\ep^* [f])_k  \rightarrow \op{L}^*f(x_k)$, as desired. But for general $q$ this is not the case, and since we assumed no knowledge of $q$ we have no way to control the extra terms that corrupt the approximation of $\op{L}$ and $\op{L}^*$. In Section \ref{sec:localkernel_unbiased} we will show how to modify the construction by Berry and Sauer so that this unsatisfactory situation can be remedied.

We comment briefly on the origin of the extra drift terms in Theorem \ref{thm:localkernels_q}. As $m\rightarrow \infty$, we have
\[
K_\ep^{A,b} [f]_k \rightarrow \int_{\op{M}} k^{A,b}_\ep(x_k, y) q(y) f(y) dy
\]
in the Monte Carlo sense, with a $\mathcal{O}(m^{-1/2})$ variance term whose exact form is reported in \cite{SINGER2006}. On the other hand, Berry and Sauer \cite{BeSa16} show the integral operator expansion
\begin{align*}
\op{G}_\ep f(x) &\; \equiv \ep^{-d/2}\int_M k^{A,b}_\ep (x,y)f(y) dy \\
&\; = m(x)f(x) + \ep\left[\omega(x)f(x) + m(x)\op{L}f(x)\right] + \mathcal{O}(\ep^{3/2}),
\end{align*}
where $m(x)$ and $\omega(x)$ are scalar functions depending on the kernel $k^{A,b}_\ep$ and the induced metric on $\op{M}$. Thus, up to the normalization factor $\ep^{-d/2}$ which is eventually cancelled out, $K_\ep^{A,b}[f]$ approximates $\op{G}_\ep(fq)$ and not $\op{G}_\ep f$. Details can be found in Appendix \ref{sec:preliminaries}.


\section{Local kernel diffusion maps}\label{sec:localkernel_unbiased}

With the same notation as in the previous section, we propose the \emph{local kernel diffusion maps} (LKDmap) algorithm which constructs two $m\times m$ matrices $L_\ep$ and $L^*_\ep$ on the data $\set{D}^{(m)}$ that approximate the differential operators $\op{L}$ and $\op{L^*}$ in \eqref{eq:backward_gen} and \eqref{eq:forward_gen}. The algorithm works with two different kernels: The local kernel $k^{A,b}_\ep$ from \eqref{eq:localkernel} and the isotropic kernel $k_\ep$ from \eqref{eq:kernel}, which is equal to $k^{A,b}_\ep$ with $A = I$ and $b=0$. The second kernel is used in order to construct a kernel density estimate $q_\ep$ of the unknown sampling density $q$ in order to remove the resulting unwanted drift terms.

As in TMDmap, the idea is to normalize the kernel matrix $K_\ep^{A,b}$ twice: The right normalization by the diagonal matrix $D_{\tilde\ep}$ in step (ii) removes unwanted drift terms produced by local variations in the sampling density $q$. This is done by means of a kernel density estimator $q_{\tilde\ep}$ constructed from the isotropic kernel $k_{\tilde \ep}$. The left normalization in step (iii) ensures that probability is conserved.

\begin{center}
\colorbox{boxback}{
  \begin{minipage}[t]{\textwidth}
\begin{center}
  \begin{minipage}[t]{0.95\textwidth}
\vspace{0.5cm}
\center\textbf{Construction of LKDmap}
\begin{enumerate}[(i)]
\item Construct the $m\times m$ kernel matrices $K^{A,b}_\ep$ and $K_{\tilde\ep}$ with components $(K^{A,b}_\ep)_{ij} = k^{A,b}_\ep(x_i, x_j)$ and $(K_{\tilde\ep})_{ij} = k_{\tilde \ep}(x_i, x_j)$, and the kernel density estimate
\[
q_{\tilde\ep}(x_i) = \sum_{j=1}^m (K_{\tilde\ep})_{ij}.
\]
\item Form the diagonal matrix $D_{\tilde\ep}$ with components $(D_{\tilde\ep})_{ii} = q^{-1}_{\tilde\ep}(x_i)$ and set
\[
\tilde K^{A,b}_{\ep,\tilde\ep} = K^{A,b}_\ep D_{\tilde\ep},\quad (\tilde K^{A,b}_{\ep,\tilde \ep})^* = (K^{A,b}_\ep)^T D_{\tilde\ep}.
\]
\item Let $\tilde D_{\ep,\tilde \ep}$ be the diagonal matrix of row sums of $\tilde K^{A,b}_{\ep,\tilde \ep}$, that is,
\[
(\tilde D_{\ep,\tilde \ep})_{ii} = (\tilde K^{A,b}_{\ep,\tilde \ep} [\mathbf{1}])_i = \sum_{j=1}^m (\tilde K^{A,b}_{\ep,\tilde \ep})_{ij}.
\]
\item Build the LKDmap matrices
\begin{subequations}
\begin{align}
L_{\ep,\tilde \ep} &\;= \ep^{-1}\left( \tilde D_{\ep,\tilde \ep}^{-1} \tilde K^{A,b}_{\ep,\tilde \ep} - I\right), \label{eq:LKDmap}\\ 
L^*_{\ep,\tilde \ep} &\;= \ep^{-1}\left((\tilde K^{A,b}_{\ep,\tilde \ep})^* \tilde D_{\ep,\tilde \ep}^{-1} - I\right) \label{eq:LKDmap_star}.
\end{align}
\end{subequations}
\end{enumerate}
  \end{minipage}
\vspace{0.5cm}
\end{center}
 \end{minipage}
}
\end{center}

We note that $(D_{\tilde \ep}L_{\ep,\tilde\ep})^T = D_{\tilde\ep}L^*_{\ep,\tilde \ep}$, which naturally extends the relation $L_\ep^T = L_\ep^*$ from section \ref{sec:anisotropic} to the non-uniform sampling case. Since LKDmap uses two kernels, there are also two associated scale parameters $\ep$ and $\tilde \ep$. The interpretation of the two parameters is rather different: One can think of $\ep$ as a characteristic time scale for the dynamics \eqref{eq:SDE} associated with $A$ and $b$; this is apparent in \eqref{eq:localkernel} as the term $x+\ep b(x)$ mimics an explicit Euler step. On the other hand, $\tilde \ep$ is best thought of as the width of the convolution kernel that is used to construct the kernel density estimate $q_{\tilde \ep}$. It is therefore not always appropriate (and not necessary) to choose $\tilde \ep = \ep$, and we will keep both parameters separate in our analysis. Limiting results are obtained for $\ep\rightarrow 0$ and $\tilde \ep \rightarrow 0$.

The following theorem establishes the relation between the matrices $L_\ep$ and $L_\ep^*$ and the differential operators $\op{L}$ and $\op{L}^*$. The proof can be found in Appendix \ref{sec:proof_LKDmap}.

\begin{theorem}[Convergence of LKDmap]
\label{thm:local_diffmap}
Let $A(x)$ be a matrix-valued function on $\op{M}\subset \R^N$ such that each $A(x)$ is a symmetric positive definite $N\times N$ matrix, and let $b(x)$ be a vector-valued function. Let $L_{\ep,\tilde\ep}$ and $L_{\ep,\tilde\ep}^*$ be defined via \eqref{eq:LKDmap} and \eqref{eq:LKDmap_star}, respectively, and let $f : \op{M} \rightarrow \R$ be smooth. Then in the limit $m\rightarrow \infty$, $\ep \rightarrow 0$ and $\tilde\ep \rightarrow 0$ and for every $x_k \in \set{D}^{(m)}$, we obtain
\begin{align*}
(L_{\ep,\tilde \ep} [f])_k &\; \rightarrow \op{L}f(x_k),\\
(L_{\ep,\tilde \ep}^* [f])_k &\; \rightarrow \op{L}^*f(x_k),
\end{align*}
where $\op{L}$ and $\op{L}^*$ are defined in \eqref{eq:backward_gen} and \eqref{eq:forward_gen}, respectively.
\end{theorem}
The convergence in Theorem \ref{thm:local_diffmap} is of the same nature as in Theorem \ref{thm:target_diffmap} and standard diffusion maps \cite{CoLa06}. As in Theorem \ref{thm:target_diffmap}, solutions to the matrix eigenvalue problems $L_\ep [\psi] = \lambda [\psi]$ or $L_\ep^*[\psi] = \lambda[\psi]$, respectively, can be used to approximate solutions of the corresponding eigenvalue problems for $\op{L}$ and $\op{L}^*$ on the domain $\op{M} = \supp(q)$, with von Neumann boundary conditions imposed on $\partial \op{M}$ if $\op{M}$ has a boundary.

LKDmap requires the knowledge of $A(x_i)$ and $b(x_i)$ for all $x_i \in \set{D}^{(m)}$. If $A$ and $b$ are unknown but one has access to a blackbox integrator that produces realizations of a dynamical process $X_t$ with prescribed initial condition $X_0$, then one can construct estimators of $A(x)$ and $b(x)$ based on the Kramers--Moyal expansion \cite{risken1996fokker}
\begin{subequations}
\begin{align}
b(x) &\; = \lim_{\tau\rightarrow 0}\frac{1}{\tau}\mathbf{E}\left[X_\tau - X_0 | X_0 = x\right] \label{eq:Kramers_b}\\
A(x) &\; = \lim_{\tau\rightarrow 0}\frac{1}{2\tau}\mathbf{Cov}\left[X_\tau - X_0 | X_0 = x\right]\label{eq:Kramers_A}
\end{align}
\end{subequations}
by e.g.~many independent realizations of $X_t$ starting at $X_0 = x$ with short integration time $\tau$. If the dynamical process $X_t$ is an It\^o diffusion of the form \eqref{eq:SDE}, then \eqref{eq:Kramers_b} and \eqref{eq:Kramers_A} will recover the exact drift and diffusion coefficients of the dynamics \eqref{eq:SDE}. If not, then estimators based on \eqref{eq:Kramers_b} and \eqref{eq:Kramers_A} will compute the drift and diffusion coefficients of the best approximation of the unknown dynamics $X_t$ by an It\^o diffusion of the form \eqref{eq:SDE} in the sense of projected dynamics \cite{zhang2017effective}.

Numerical estimators $\hat A(x)$ and $\hat b(x)$ of $A(x)$ and $b(x)$ will be subject to statistical error, which can corrupt the matrix inverse in \eqref{eq:localkernel}. If one is working with numerical estimators, it may therefore be necessary to regularize the inverse in the computation of $k_\ep^{A,b}$, e.g.~by replacing \eqref{eq:localkernel} with the regularized version
\begin{equation}
k^{A,b}_{\ep,\eta}(x,y) = \exp\left(-\left(4\ep\right)^{-1} \left(x-y+\ep \hat b(x)\right)^{T}(\hat A(x) + \eta I)^{-1}\left(x-y+\ep \hat b(x)\right)\right).
\label{eq:localkernel_reg}
\end{equation}
The regularization parameter $\eta$ may be chosen based on estimates of the statistical error of $\hat A(x)$.

\begin{remark}
If one wishes to approximate the Kolmogorov operator $\op{L}f = \beta^{-1}\Delta f - \nabla U \cdot \nabla f$, one now has two options: The first is to use TMDmap with $\pi \propto \exp(-\beta U)$, which requires the evaluation of $U$ on the data $\set{D}^{(m)}$. The corresponding matrix $L_{\ep,\pi}$ then approximates $\beta \op{L}$ due to Theorem \ref{thm:target_diffmap}. The second option is to use LKDmap with $A = \beta^{-1}I$ and $b = -\nabla U$. This requires the evaluation of $\nabla U$ on the data $\set{D}^{(m)}$. The corresponding matrix $L_\ep$ then approximates $\op{L}$ due to Theorem \ref{thm:local_diffmap}. In our numerical examples, we do not see any significant difference in approximation quality between TMDmap and LKDmap in this situation. If evaluating $U$ is significantly cheaper than evaluating $\nabla U$, then TMDmap might be preferable.
\end{remark}

\section{Numerical examples}\label{sec:numerics}

We discuss a range of numerical examples. Matlab code that implements TMDmap and LKDmap is available online at \url{https://github.com/ZofiaTr/GeneralizedDiffusionMap}.

\subsection{Removing large time step bias}

\begin{figure}[h]
\centering
\includegraphics[width=0.48\textwidth]{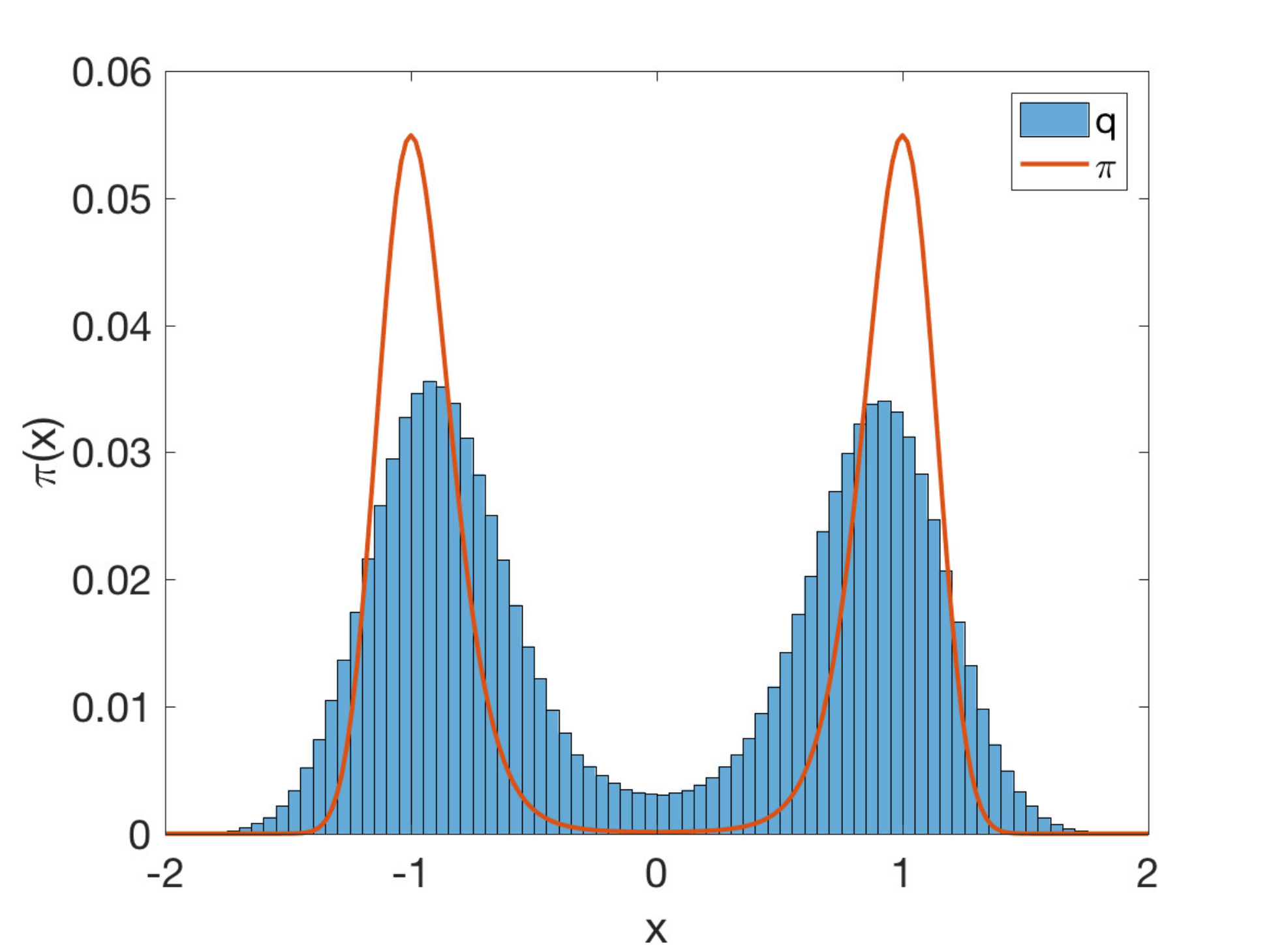}
\includegraphics[width=0.48\textwidth]{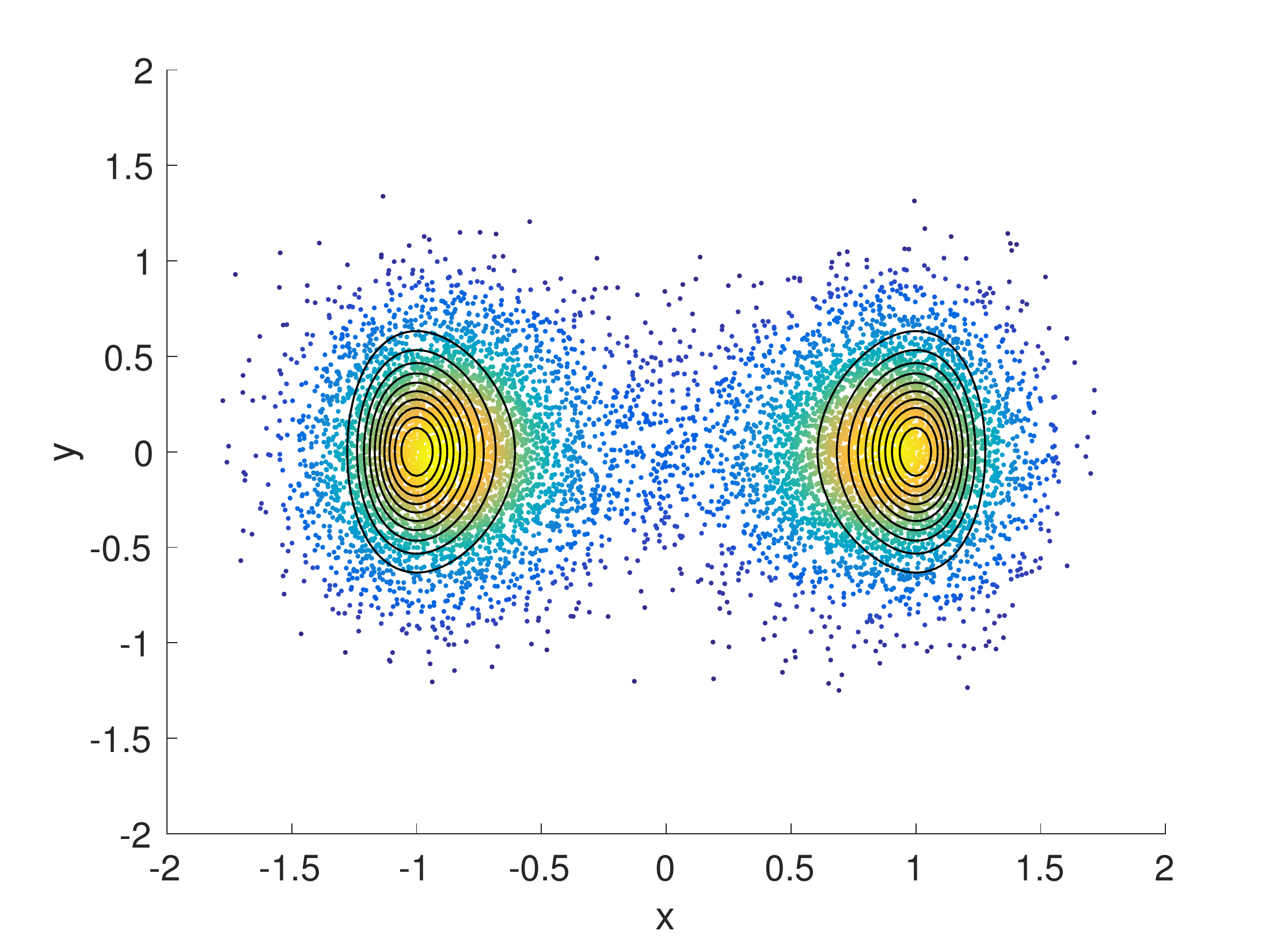}
\includegraphics[width=0.48\textwidth]{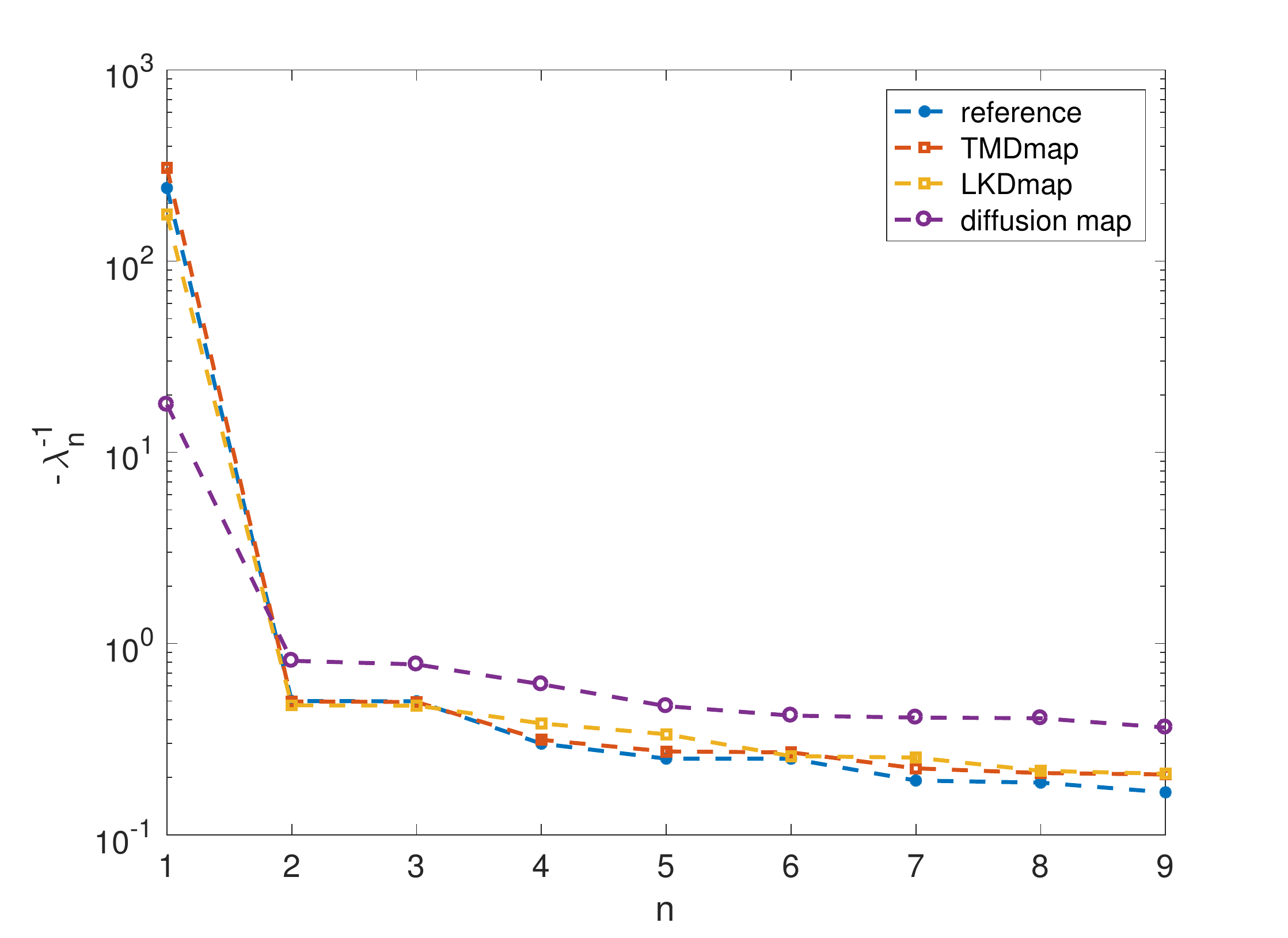}
\includegraphics[width=0.48\textwidth]{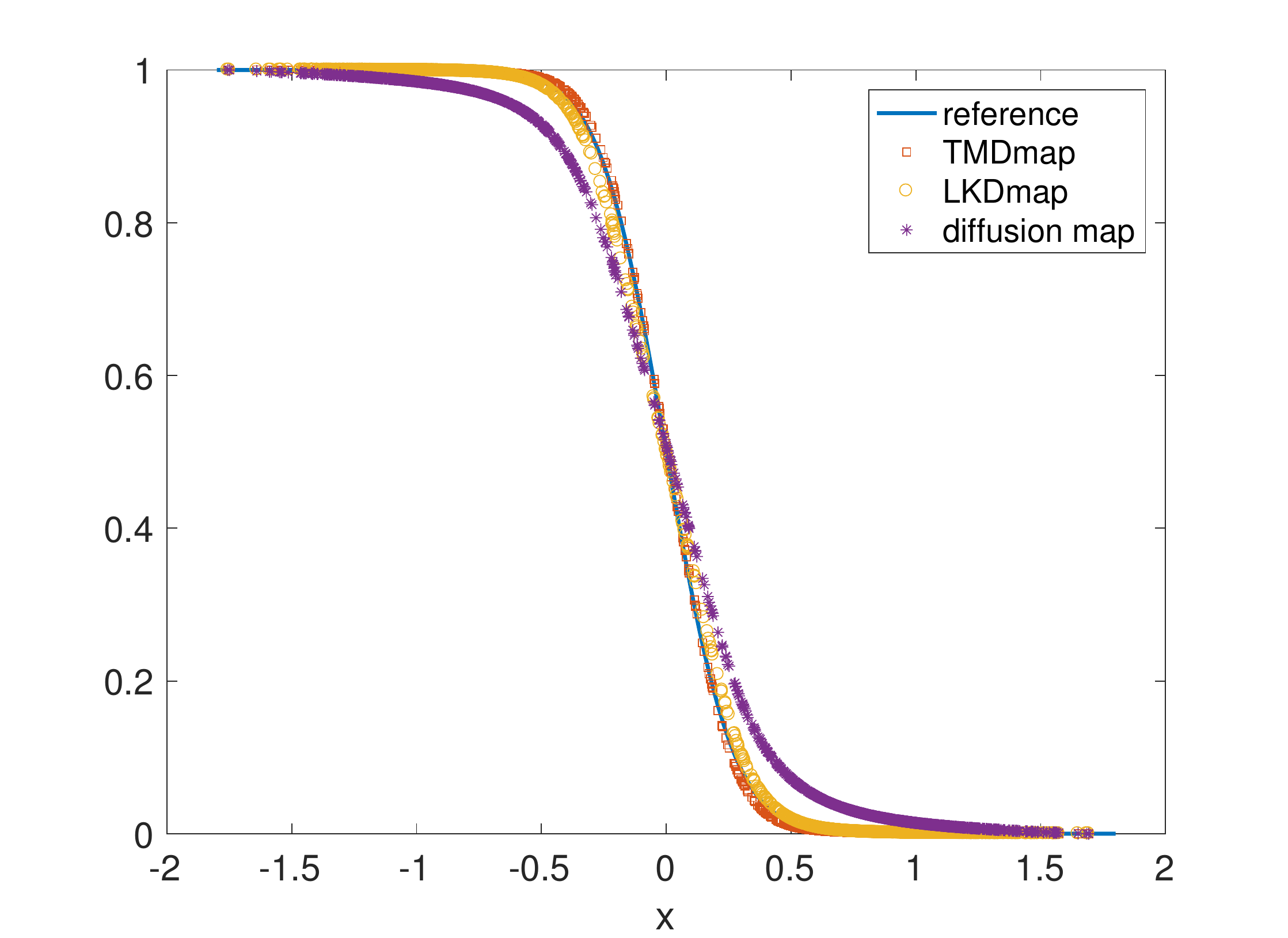}
\caption{Top left: The marginal of the target distribution $\pi \propto \exp(-\beta U)$ (red) and a normalized histogram of the sampled distribution $q$ (blue) along the $x$-axis. The sampled distribution is biased due to the large timestep. Top right: A contour plot of $\pi$ (black lines) and the sampling points colored according to $q_\ep$, showing that $q_\ep$ effectively detects the bias. Bottom: Time scales $-\lambda_n^{-1}$ (left) and the dominant eigenfunction $\psi_1(x,0)$ (right) of TMDmap, LKDmap, diffusion map, and a reference discretization. Both TMDmap and LKDmap show excellent agreement with the reference result, while standard diffusion map shows errors.}
\label{fig:largedt}
\end{figure}

In this section, we demonstrate that both TMDmap and LKDmap are capable of removing large time step bias. We study the generator \eqref{eq:backward_gen_pi} of the SDE \eqref{eq:SDE_ergodic} in dimension $N=2$ with $\pi\propto \exp(-\beta U)$ being the Boltzmann distribution with parameter $\beta = 6.0$ and the double well potential $U(x,y) = (x^2 - 1)^2 + y^2$. We generate $10^6$ samples by discretizing the SDE \eqref{eq:SDE_ergodic} with a forward Euler scheme and a time step of $\Delta t = 0.03$. This introduces a discretization bias, and as a result the sampled density $q$ is different from $\pi$ \cite{mattingly2002ergodicity}. In our case, $\Delta t$ is so large that the bias is very apparent, as can be seen in Figure \ref{fig:largedt} (top left).

We now demonstrate that TMDmap and LKDmap are capable of removing this bias. For the purpose of this computation we thin the data to $m = 10^4$ samples. Then we construct the TMDmap matrix $L_{\ep,\pi}$, which approximates the operator \eqref{eq:backward_gen_pi}. We also construct the LKDmap matrix $L_{\ep,\tilde \ep}$ with $A = \beta^{-1}I$ and $b = -\nabla U$, which approximates the same operator \eqref{eq:backward_gen_pi} up to a factor of $\beta$. Finally, we build a standard diffusion maps matrix $L_{\ep,\alpha}$ with $\alpha = 1/2$. We choose $\ep = \tilde \ep = 0.05$ in all cases.

Figure \ref{fig:largedt} compares the largest time scales $t_n = -\lambda_n^{-1}$ associated with the smallest (in magnitude) eigenvalues $\lambda_n$ of all the matrices constructed. A reference obtained with a high resolution finite element discretization \cite{latorre2011structure} is also shown. We observe that both TMDmap and LKDmap estimate the first nine time scales very well. Of particular interest is the largest time scale $t_1$ which is the characteristic jumping time between the two minima $(\pm 1, 0)$ of $U$. Standard diffusion maps on the other hand underestimates $t_1$ by one order of magnitude and overestimates the other time scales significantly. This is a direct consequence of the bias in $q$. The dominant eigenfunction $\psi_1$ associated with $\lambda_1$ is also shown in Figure \ref{fig:largedt} for all three methods. Again the FEM reference, TMDmap and LKDmap all agree very well, while diffusion maps produces significant errors while still reproducing the correct qualitative behavior.

The reason why TMDmap and LKDmap are able to remove the bias in the sampling density $q$ is that it is detectable through the kernel density estimate $q_\ep$. To illustrate this, the data points colored with $q_\ep$ overlaid with a contour plot of $\pi$ is shown in Figure \ref{fig:largedt} in the top right. It is evident that $q_\ep$ approximates $q$ well and thus effectively detects the bias introduced by~$q$. An alternative approach based on a least squares approximation of the evolution operator associated with the dynamics was considered in \cite{wu2017variational, KlKoSch16}.

\subsection{Temperature switch}

\begin{figure}[h]
\centering
\includegraphics[width=0.45\textwidth]{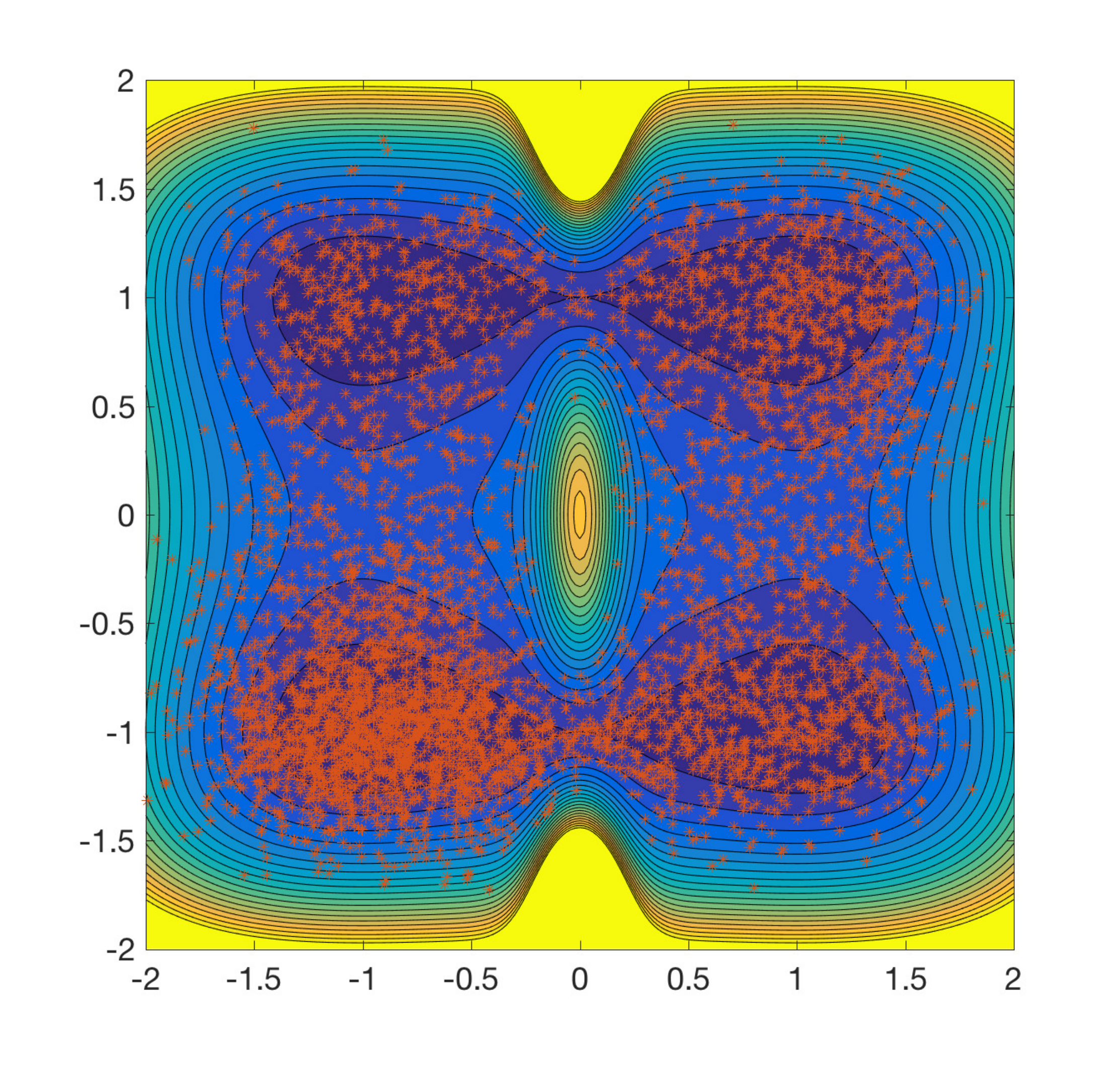}
\caption{The potential $U$ from \eqref{eq:pot_switch} and $m=5,000$ sampling points sampled from a realization of the dynamics \eqref{eq:SDE} with $b = -\nabla U$, $\sigma^2 = 2\beta_s^{-1}$ with $\beta_s = 1.0$ and initial condition $X_0 = (-1,-1)$.}
\label{fig:switch1}
\end{figure}

\begin{figure}[h]
\centering
\includegraphics[width=0.9\textwidth]{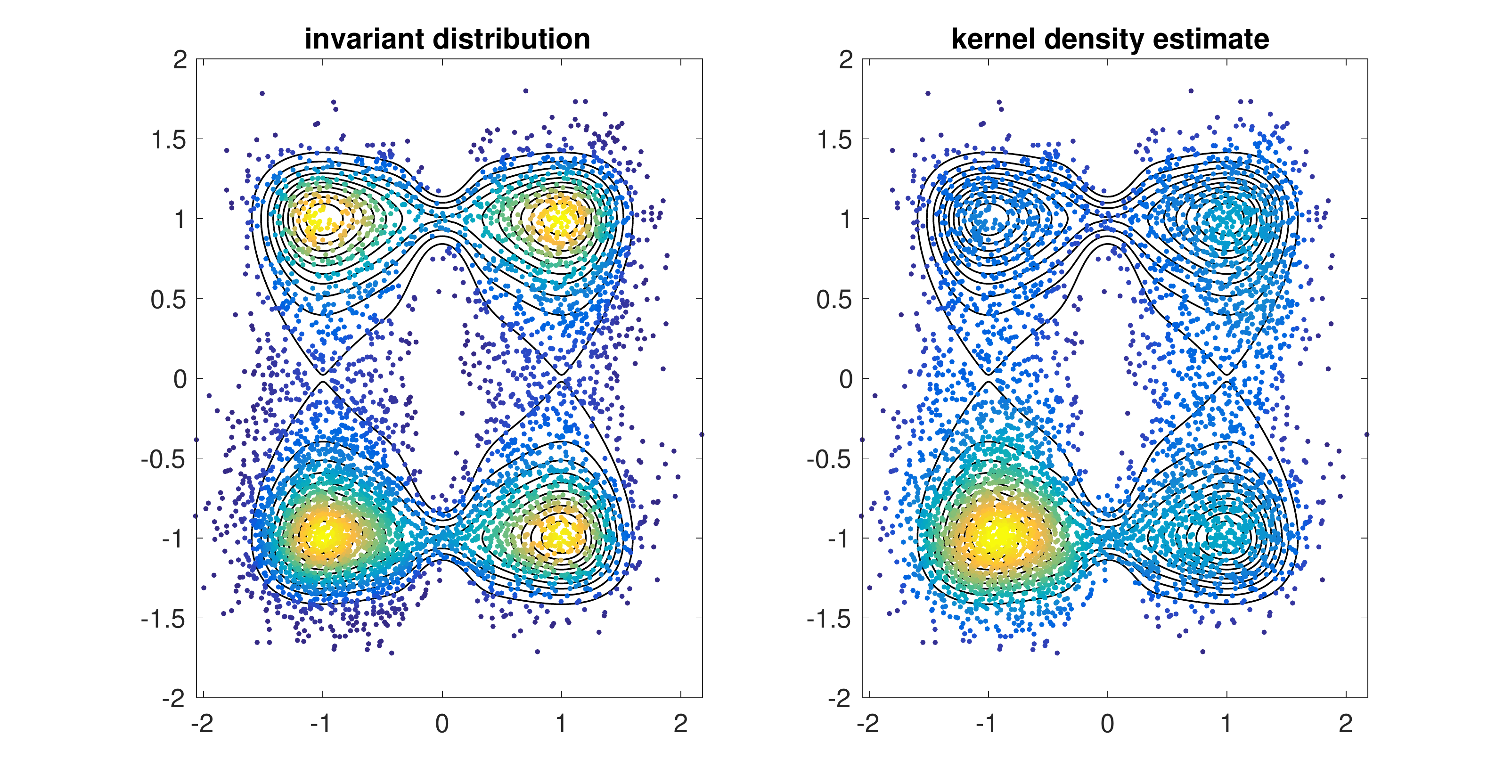}
\caption{Left: The left zero eigenvector $\phi_0$ of the TMDmap matrix $L_{\ep,\pi}$ times the kernel density estimate $q_\ep$ (colored points) is a good estimator of the invariant density $\pi$ for $\beta = 2.0$ (contour lines). The relative error in $l_1$ norm is $7.2\%$. Right: The same points colored according to $q_\ep$ show that the data is not equilibrated and sampled from the higher temperature $\beta_s = 1.0$. The relative error in $l_1$ norm between $q_e$ and $\pi$ is $43 \%$.}
\label{fig:switch2}
\end{figure}

We demonstrate that TMDmap and LKDmap are capable of effectively approximating the generator \eqref{eq:backward_gen_pi} with $\pi \propto \exp(-\beta U)$ even in cases where samples are generated using a different temperature $\beta_s \neq \beta$, and are not converged. We study the potential
\begin{equation}
\label{eq:pot_switch}
U(x,y) = h_x (x^2-1)^2 + (h_y + a(x,\delta))(y^2 - 1)^2
\end{equation}
in $N=2$ dimensions with
\[
a(x,\delta) = \frac{1}{5}\left( 1 + 5 \exp(-(x-x_0)^2/\delta)\right)^2
\]
and the parameters $h_x = 0.5$, $h_y = 1.0$, $x_0 = 0$, and $\delta = 1/20$. The potential is shown in Figure~\ref{fig:switch1} along with $m = 5,000$ sampled points from a realization of the dynamics \eqref{eq:SDE} with $b = -\nabla U$, $\sigma = \sqrt{2\beta_s^{-1}}$, time step $\Delta t = 0.02$, and initial condition $X_0 = (-1,-1)$. This dynamics samples the density $q \propto \exp(-\beta_s U)$ at the high temperature $\beta_s^{-1} = 1.0$. Figure \ref{fig:switch1} shows that the samples have explored most of the state space, but are clearly not equilibrated. Figure \ref{fig:switch2} on the right shows the kernel density estimate $q_\ep$, which is highest in the bottom left minimum where the initial condition is located.

Additionally, the function $U$ is designed so that the dynamics \eqref{eq:SDE} with $b = -\nabla U$ and $\sigma = \sqrt{2\beta^{-1}}$ shows qualitatively different behavior for different values of $\beta$. $U$ has four minima at $(\pm 1, \pm 1)$. The energy barrier in the $y$ direction is higher then in the $x$ direction, but the additional coupling term $a(x,\delta)$ creates an entropic barrier in the $x$ direction. For small $\beta$, the entropic barrier in $x$ direction is harder to cross then the energy barrier in the $y$ direction, and we expect the dominant eigenfunction $\psi_1$ of $\op{L}$ to parameterize the slow $x$  variable. For large $\beta$, the energy barrier is harder to cross, and we expect $\psi_1$ to parameterize the $y$ coordinate. Thus, $\beta$ acts as a switch which completely changes the behavior of the dominant eigenfunctions.

We construct the TMDmap matrix $L_{\ep,\pi}$ for $\ep = 0.2$ and compute the left zero eigenvector $\phi_0$. The left of Figure \ref{fig:switch2} shows the estimate of the target distribution $\pi = \exp(-\beta U)$ according to \eqref{eq:pi_hat} for $\beta = 2.0$. The truth is shown as contour lines. The approximation seems reasonable.

\begin{figure}[ht]
\centering
\includegraphics[width=0.48\textwidth]{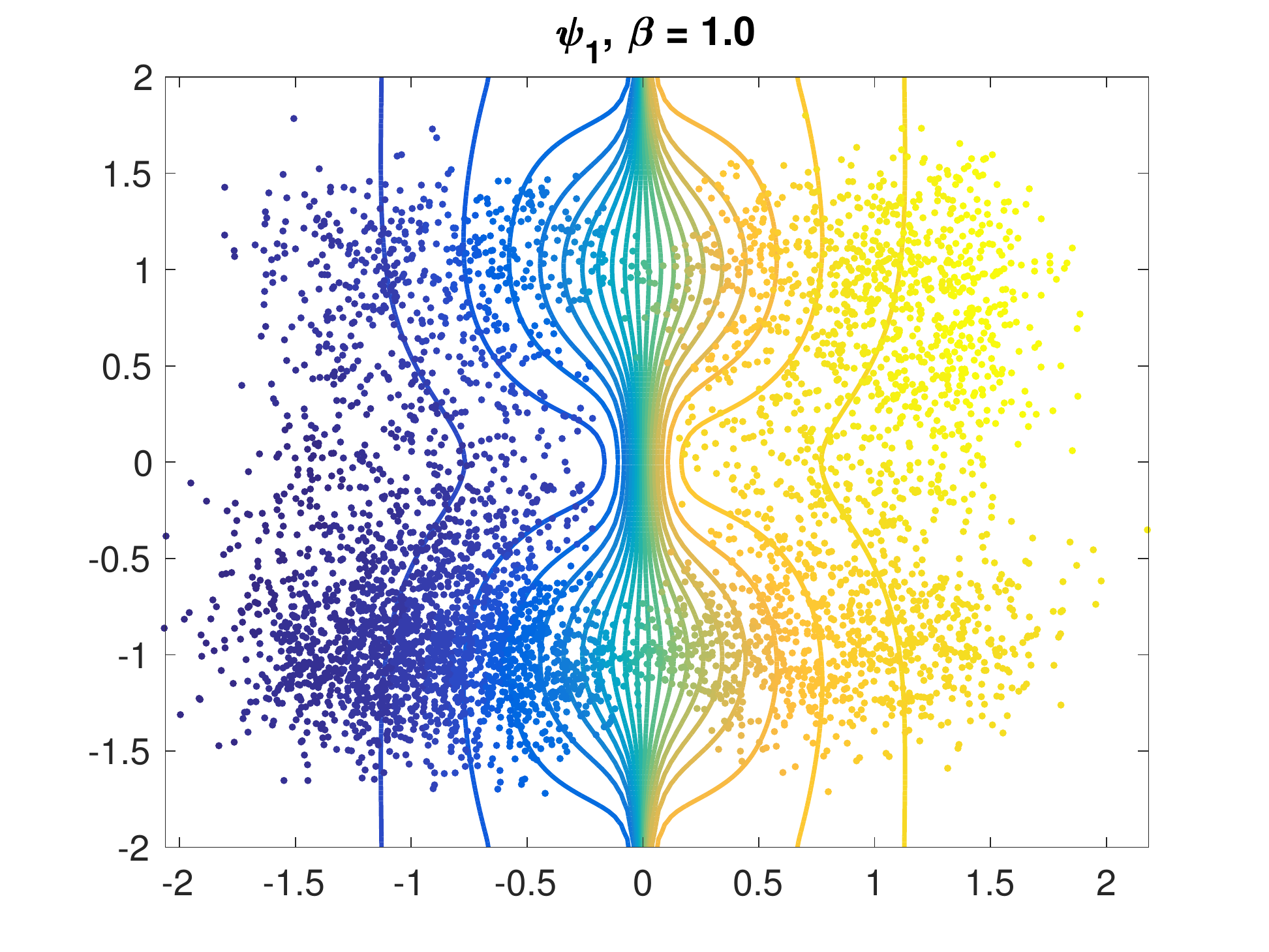}
\includegraphics[width=0.48\textwidth]{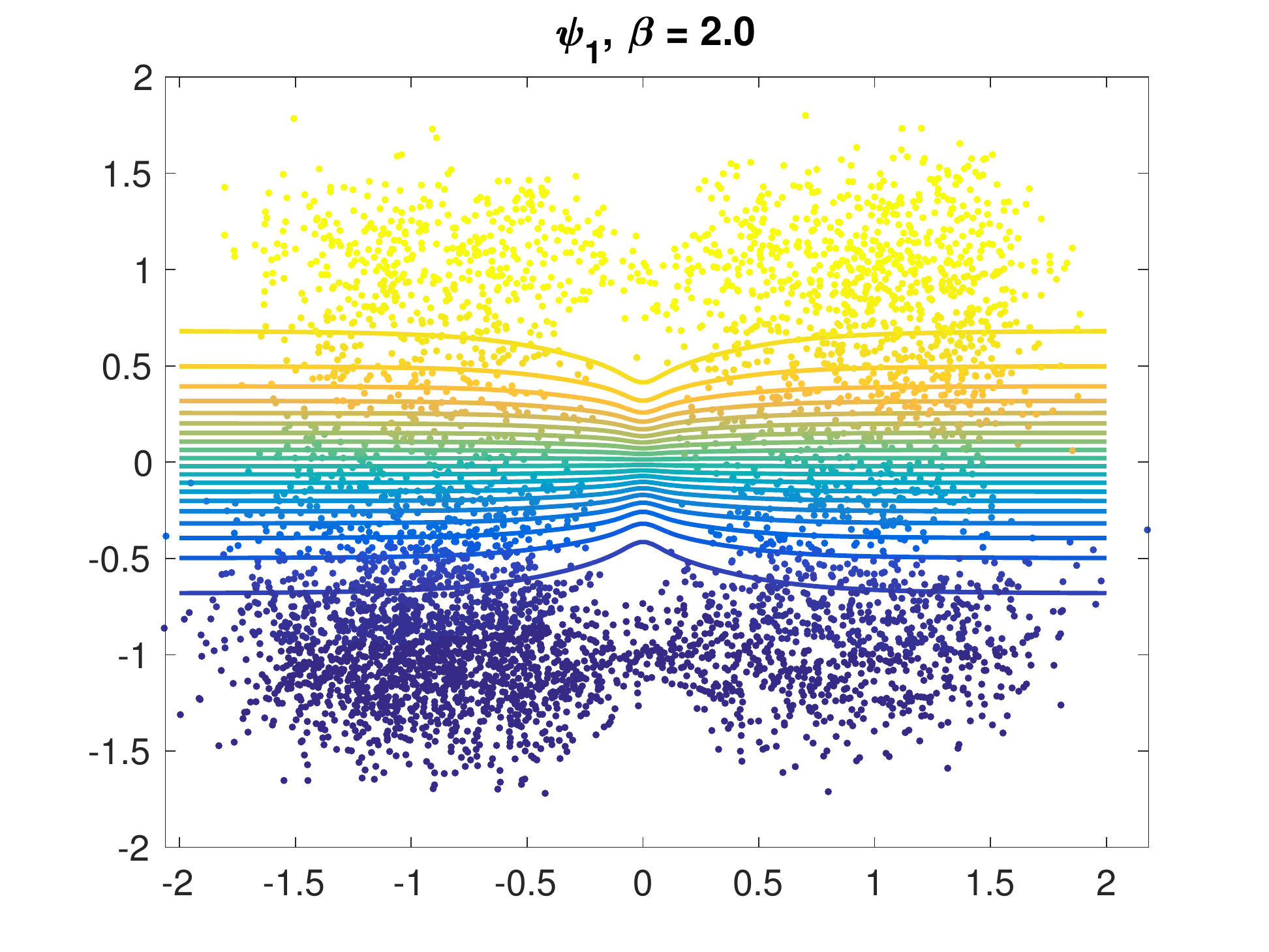}
\caption{Left: The dominant eigenvector of the TMDmap matrix $L_{\ep,\pi}$ for $\beta = 1.0$. The contour lines represent the FEM solution. Right: The dominant eigenvector of $L_{\ep,\pi}$ for $\beta = 2.0$. The relative error in $l_2$ norm, computed by cubic interpolation of the FEM solution, is $4.8\%$ and $5.5\%$ respectively.}
\label{fig:switch3}
\end{figure}

\begin{figure}[ht]
\centering
\includegraphics[width=0.98\textwidth]{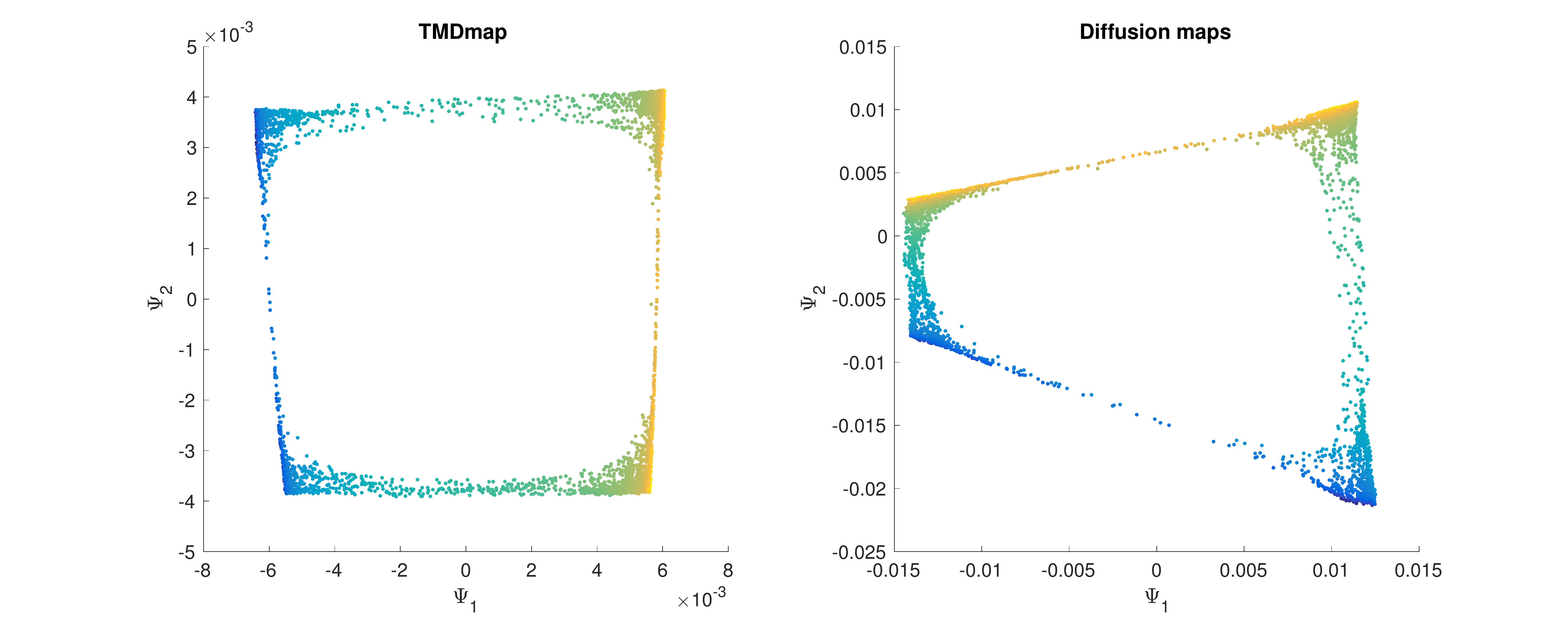}
\caption{Different embeddings for the temperature switch example \eqref{eq:pot_switch}. Left: TMDmap embedding with $\beta = 2.0$. Right: Diffusion map embedding. Data points are colored according to their location on the $y$ axis.}
\label{fig:switch_embedding}
\end{figure}

The dominant eigenvector $\psi_1$ of $L_{\ep,\pi}$ is shown in Figure \ref{fig:switch3} for $\beta = \beta_s = 1.0$ and a lower temperature $\beta = 2.0$. A reference was also computed with a FEM method \cite{latorre2011structure} on the domain $[-2,2]\times[-2,2]$ with von Neumann boundary conditions. Figure \ref{fig:switch3} shows that as expected, for the high temperature the $x$ coordinate is the slow coordinate, and for the low temperature the $y$ coordinate is the slow coordinate. TMDmap captures this in both cases even though only samples at $\beta_s = 1.0$ have been used. Figure \ref{fig:switch_embedding} compares the embeddings $\Psi(x) = (e^{\lambda_1} \psi_1(x), e^{\lambda_2}\psi_2(x))$ produced by the first two eigenfunctions of TMDmap and standard diffusion maps for $\beta = 2.0$. Data points are colored according to their location on the $y$ axis from negative (blue) to positive (yellow); note that the $y$ coordinate is the slow variable in this case. Both embeddings reproduce the circular topology of the data with the four minima of $U$ in the corners connected by narrow pathways. However, while TMDmap does recognize the $y$ coordinate as the slow coordinate ($y$ and $\psi_1$ are aligned in Figure \ref{fig:switch_embedding} on the left), diffusion maps does not. Diffusion maps also produces a more distorted embedding due to fact that the four minima of $U$ are not evenly sampled. The same computations for LKDmap with $\ep = \tilde \ep = 0.2$ give identical results (not shown).


\subsection{Effective dynamics for the butane molecule}\label{sec:numerics_c}

We now consider a more complex problem, namely the $n$-butane molecule in explicit solvent thermostatted to $300\,\textup{K}$. This molecule consists of 14 atoms and thus has a 42 dimensional configuration space. The dimension of the full system, including the solvent, is much higher. Our data $\set{D}^{(m)}$ consists of $m =  10,000 $ data points subsampled from a $100\,\textup{ns}$ trajectory computed with Amber~\cite{Amber15}. In order to remove the influence of translations and rotations of the butane molecule, we align the atoms by minimizing the root-mean-square deviation (RMSD) of the four carbon atoms.

We consider the projection onto the 12-dimensional subspace $\mathbb{X}$ of the positions of the four carbon atoms. In this subspace, the dynamics is not an It\^o diffusion of the form \eqref{eq:SDE}. However, it is known that the slow dynamics of butane is parameterized by the dihedral angle $\vartheta$, which is a nonlinear function of the coordinates of the four carbon atoms. Thus, the slow dynamics can be parametrized by a nonlinear coordinate function on~$\mathbb{X}$. It was shown in \cite{zhang2017effective} that in this situation, an effective It\^o diffusion on the slow coordinate space can be a good approximation of the slow dynamics. Let $Z_t$ be the projection of the full 42-dimensional state vector $X_t$ onto $\mathbb{X}$. The effective It\^o diffusion is of the form
\begin{equation}
dZ_t = b(Z_t)dt + \sigma(Z_t) dW_t
\label{eq:SDE_effective}
\end{equation}
with effective drift and diffusion given by analogs of the Kramers--Moyal expansion \eqref{eq:Kramers_b} and \eqref{eq:Kramers_A}
\begin{align*}
b(z) &\; = \lim_{\tau\rightarrow 0} \frac{1}{\tau}\mathbf{E}[ Z_\tau - Z_0 | Z_0 = z], \\
A(z) &\; = \lim_{\tau\rightarrow 0} \frac{1}{2\tau}\mathbf{Cov}[ Z_\tau - Z_0 | Z_0 = z],
\end{align*}
where the conditional mean and covariance are taken with respect to the marginal equilibrium density of the dynamics for a fixed value $z$ of the 12-dimensional position vector of the 4 carbon atoms. We want to approximate the generator of the SDE \eqref{eq:SDE_effective} with LKDmap. The (position dependent) effective drift and diffusion are not known and have to be estimated at the query points. This may be done by launching very short MD simulations from each query point. For simplicity, we just use the long trajectory by considering all the points that are within a radius of $r = 0.01\,\textup{nm}$ of each query point and following where these points get mapped to after a time $\tau = 0.05\,\textup{ps}$ (we cannot consider the limit $\tau \rightarrow 0$ in practice). This allows us to use $\approx 10^4 $ points to compute $b(z_i)$ and $A(z_i)$ for each query point $z_i$.

The presence of statistical errors in the estimation of $A(z_i)$ make it necessary to regularize the matrix inverse, thus we replace \eqref{eq:localkernel} by \eqref{eq:localkernel_reg}. As regularization parameter we choose $\eta = 10^{-2}$, which is of the same magnitude as the expected statistical errors in estimating $A(z_i)$. We choose $\ep = 0.05\,\textup{ps}$, which is consistent with the choice of $\tau$, and $\tilde\ep = 0.01\,\textup{nm}^2$ (note that $\ep$ naturally has units of time, while $\tilde\ep$ has units of position squared).

\begin{figure}[htb]
    \centering
    \includegraphics[width=0.95\textwidth]{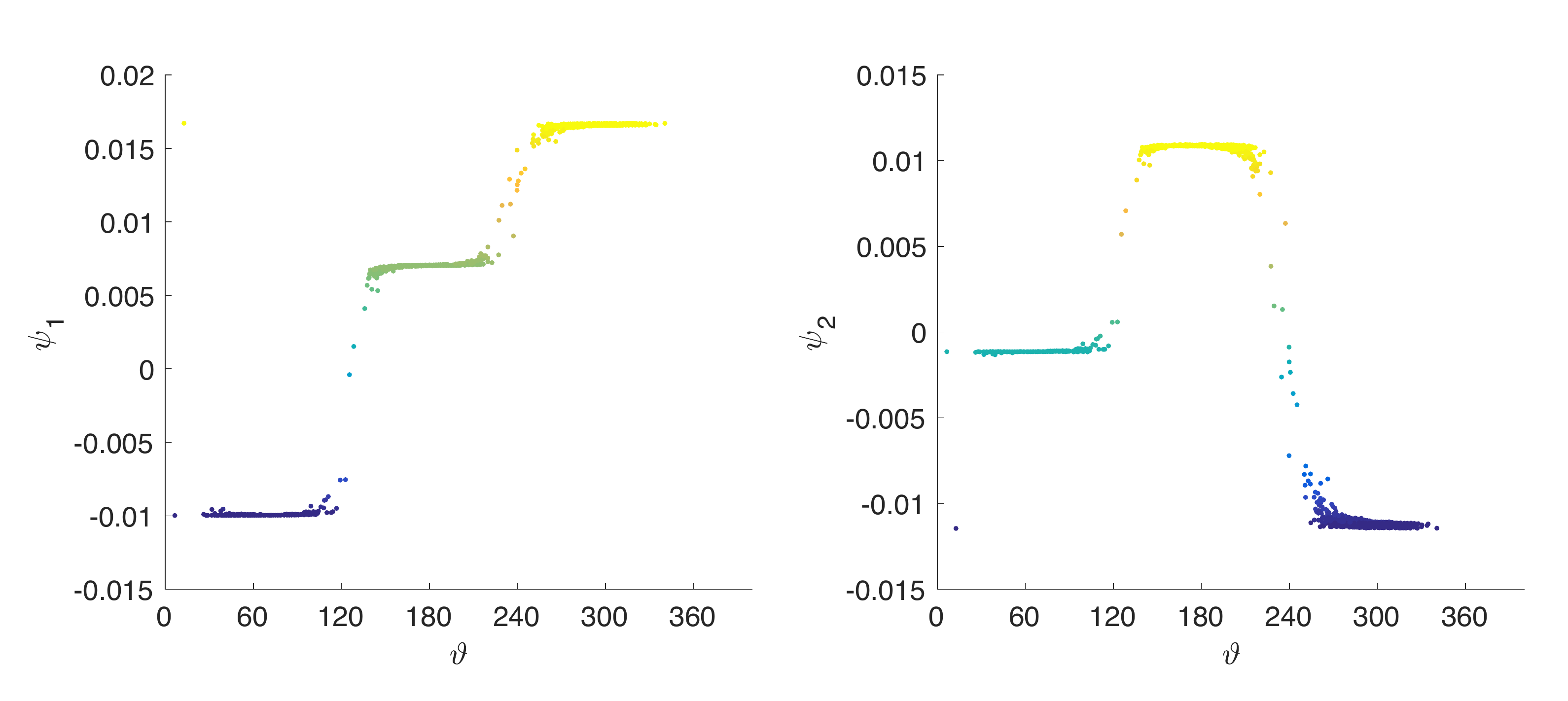}
    \caption{Eigenfunctions $\psi_1$ and $\psi_2$ of the LKDmap matrix $L_{\ep,\tilde \ep}$ associated with the butane molecule. The data points $ x_i $ were extracted from one long trajectory. The eigenfunctions are plotted as a function of the dihedral angle $ \vartheta $, which is the known slow variable.}
    \label{fig:Butane TO}
\end{figure}

\begin{figure}[htb]
    \centering
    \includegraphics[width=0.32\textwidth]{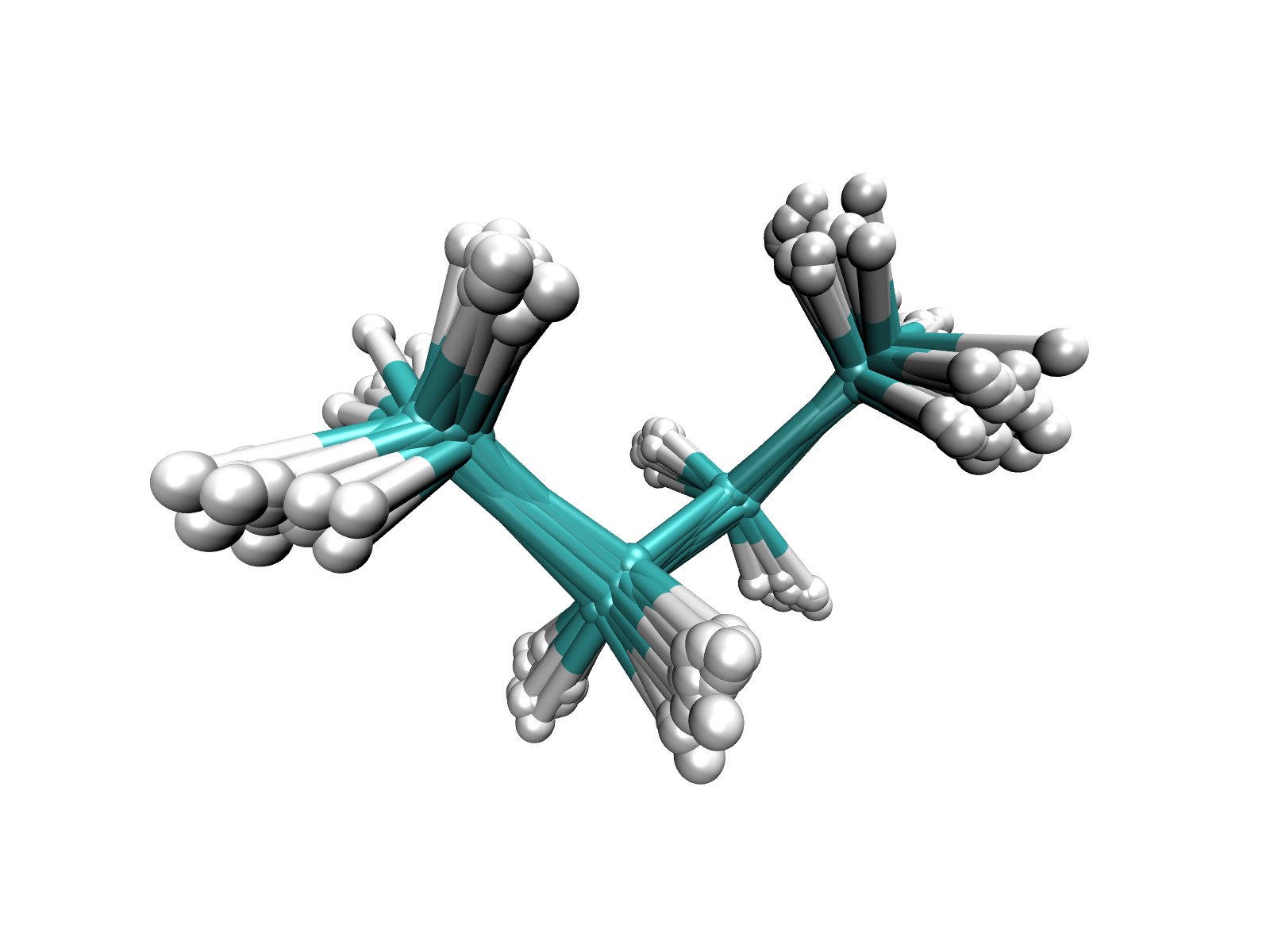}
    \includegraphics[width=0.32\textwidth]{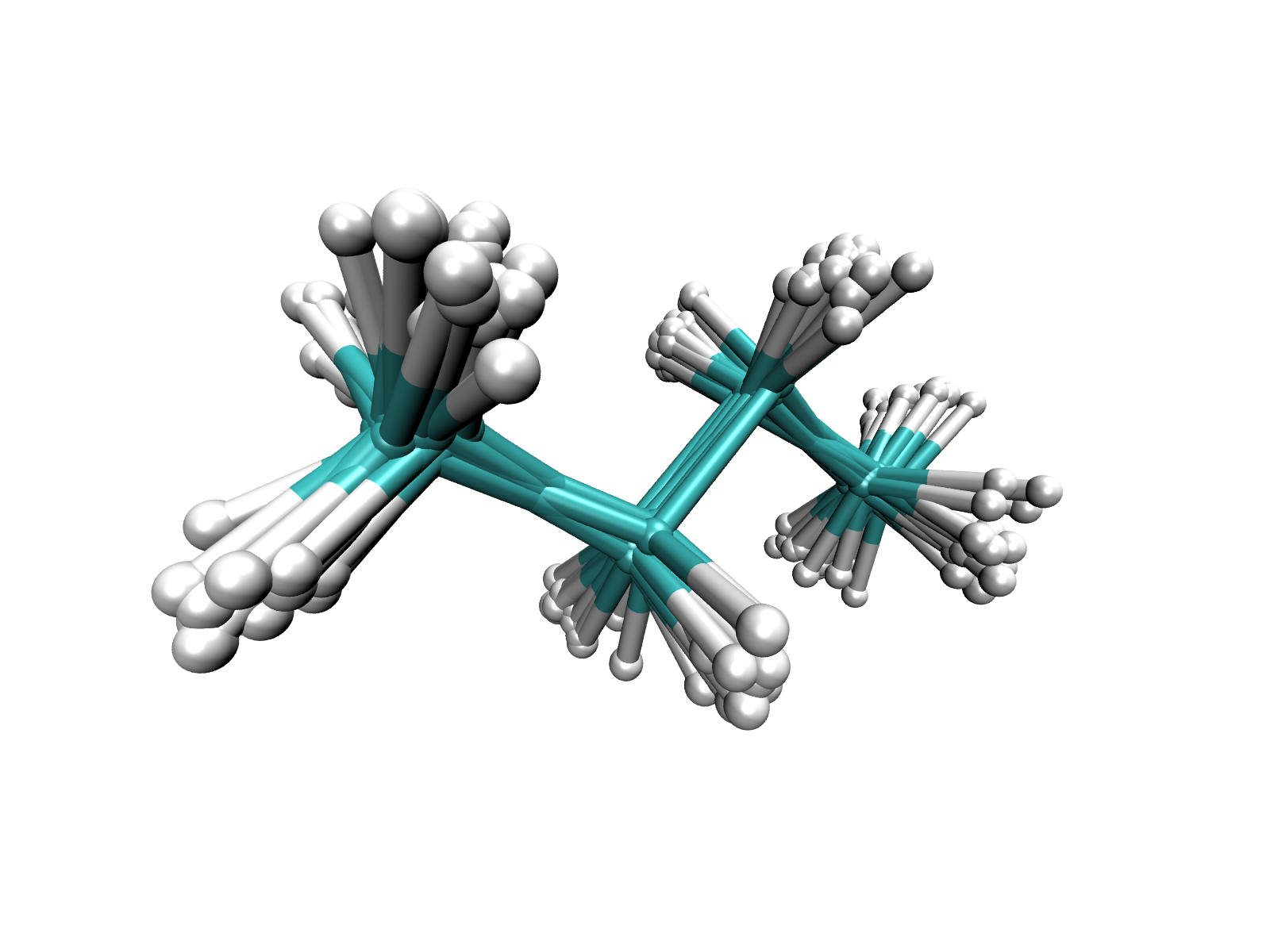}
    \includegraphics[width=0.32\textwidth]{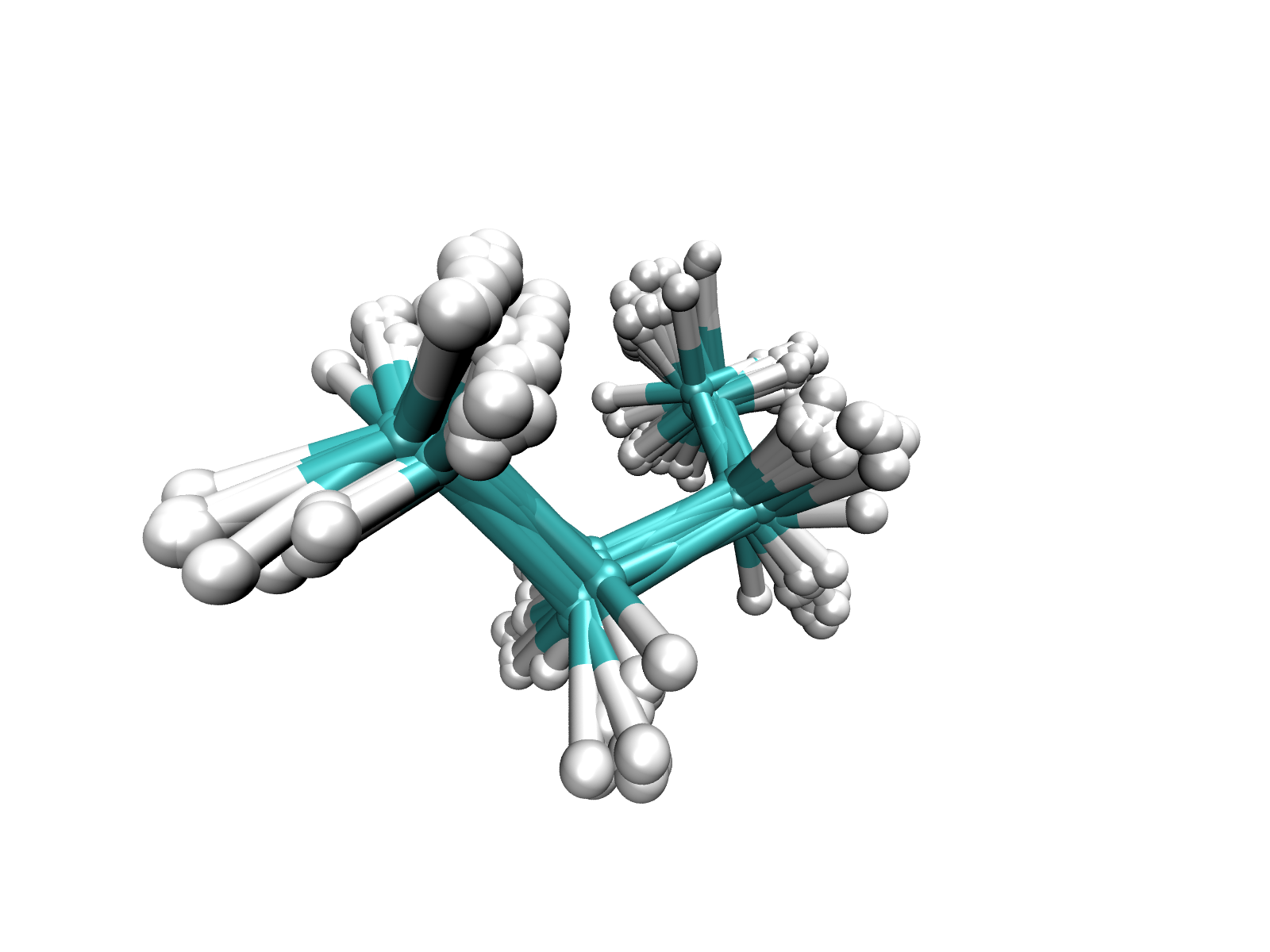}
    \caption{The configurations of the Butane Molecule associated with the three metastable states around $ \vartheta = 60^\circ $, $ \vartheta = 180^\circ $, and $ \vartheta \approx 300^\circ $.}
    \label{fig:Butane configurations}
\end{figure}

The eigenfunctions $\psi_1$ and $\psi_2$ of the LKDmap matrix $L_{\ep,\tilde \ep}$ are shown in Figure~\ref{fig:Butane TO}. We plot the values of the two eigenfunctions evaluated at all data points $ x_i $ as a function of the dihedral angle~$ \vartheta $ in order to test if LKDmap is able to discover the slow variable $\vartheta$ which is hidden in the data. We observe that the eigenfunctions $\psi_1$ and $\psi_2$ are parametrized by~$\vartheta$ (cf.~\cite{Bittracher2017} for why this is the case), and plateaus in $\psi_1$ and $\psi_2$ clearly show the expected three metastable sets around $ \vartheta = 60^\circ $, $ \vartheta = 180^\circ $, and $ \vartheta \approx 300^\circ $, corresponding to the anti and Gauche conformations. The configurations in full space corresponding to the identified metastable sets are shown in Figure~\ref{fig:Butane configurations}. Thus LKDmap is able to uncover the hidden slow variable $\vartheta$, which is nonlinear function of the given time series~$Z_t$. 

These results agree with an analysis based on Extended Dynamic Mode Decomposition~\cite{KlKoSch16}. There is a clear spectral gap between $\lambda_2$ and $\lambda_3$.  With a Markov State Model analysis with the pyemma software package~\cite{scherer2015pyemma}, we find the implied time scales $\tilde t_1 = 47\,\textup{ps}$ and $\tilde t_2 = 37\,\textup{ps}$ while the computed eigenvalues $\lambda_1$ and $\lambda_2$ suggest $t_1 = 74\,\textup{ps}$ and $t_2 = 36\,\textup{ps}$. We note that $\lambda_{1/2}$ are more sensitive to changes in the parameters $\ep, \tilde \ep$ and $\eta$ then $\psi_{1/2}$. Moreover, \cite{zhang2017effective} does not guarantee that the $\lambda_i$ are well approximated by the effective dynamics \eqref{eq:SDE_effective} even if the $\psi_i$ are, this is due to the unboundedness of the generator $\op{L}$.

\subsection{Flow structures in the wake of a vehicle}\label{sec:numerics_d}

Suppose we are given positions $x_i$ and velocities $v(x_i)$ of a number of particles. LKDmap with $b(x_i) = v(x_i)$ and $A = \beta^{-1} I$ approximates the generator \eqref{eq:backward_gen} which becomes
\[
\op{L}=v\cdot \nabla +\beta^{-1} \Delta.
\]
In the limit $\beta^{-1}\rightarrow 0$, $\op{L}$ becomes the generator of the ODE corresponding to the velocity field~$v$. For finite $\beta^{-1}$, $\op{L}$ is the generator of the corresponding SDE with added isotropic noise. If the noise is small, then the dynamics given by $\op{L}$ will move along flow lines of the velocity field much faster then orthogonal to flow lines. This enables us to find almost invariant structures by studying dominant eigenfunctions of~$\op{L}$, without requiring expensive trajectory integration~\cite{FrJuKo13}. Alternative methods can be found in \cite{berkooz1993proper, graftieaux2001combining}.

\begin{figure}[ht]
 \label{fig:unbias}
    \centering
    \begin{subfigure}[b]{0.50\textwidth}
        \includegraphics[width=\textwidth]{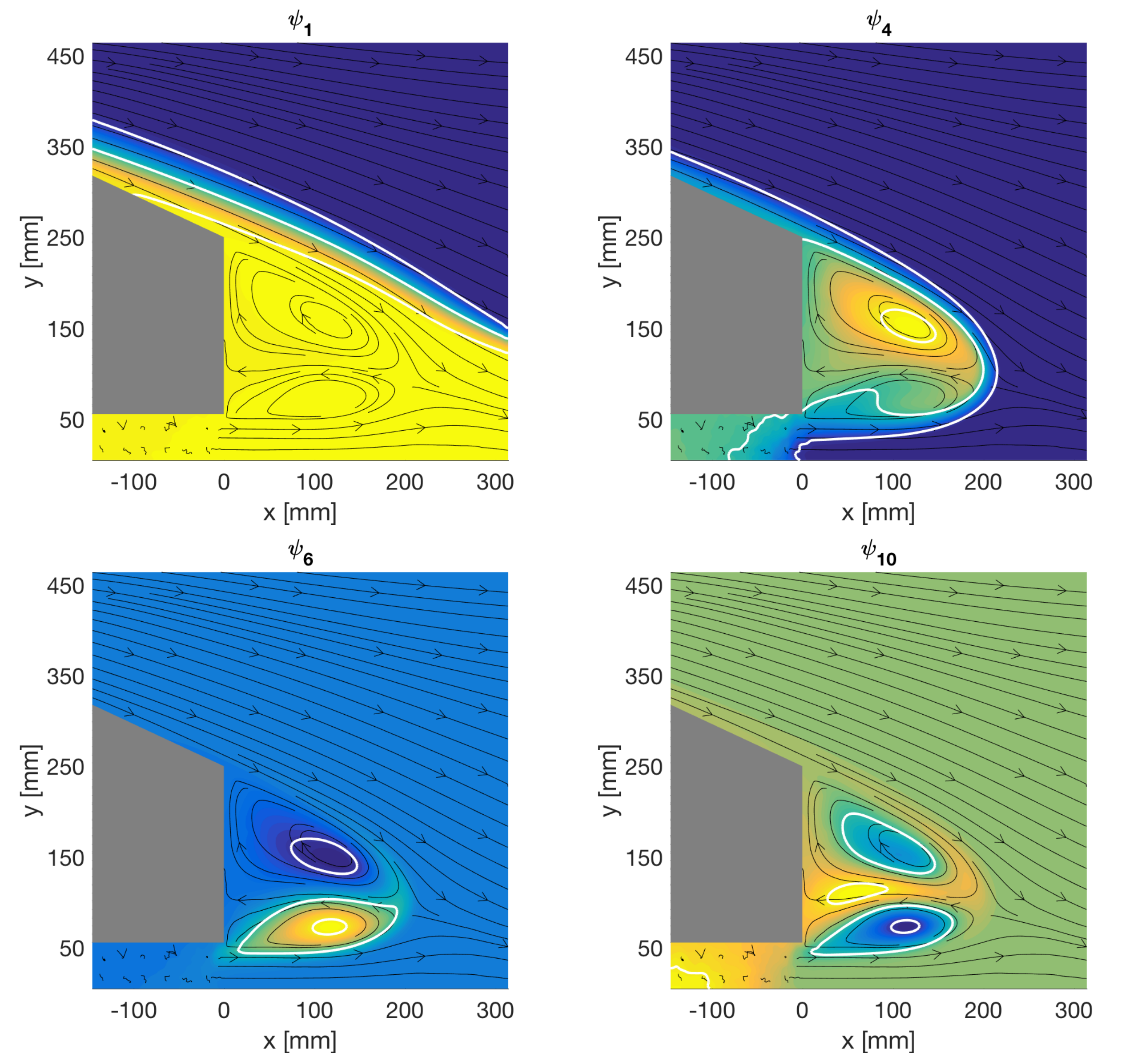}
       \caption{Data colored w.r.t indicated eigenvectors of LKDmap. Ahmed body is grey. White lines indicate the contours $[-0.9,0,0.9]$.}
        \label{fig:chosenEv_side}
    \end{subfigure}
    ~ 
    \begin{subfigure}[b]{0.46\textwidth}
        \includegraphics[width=\textwidth]{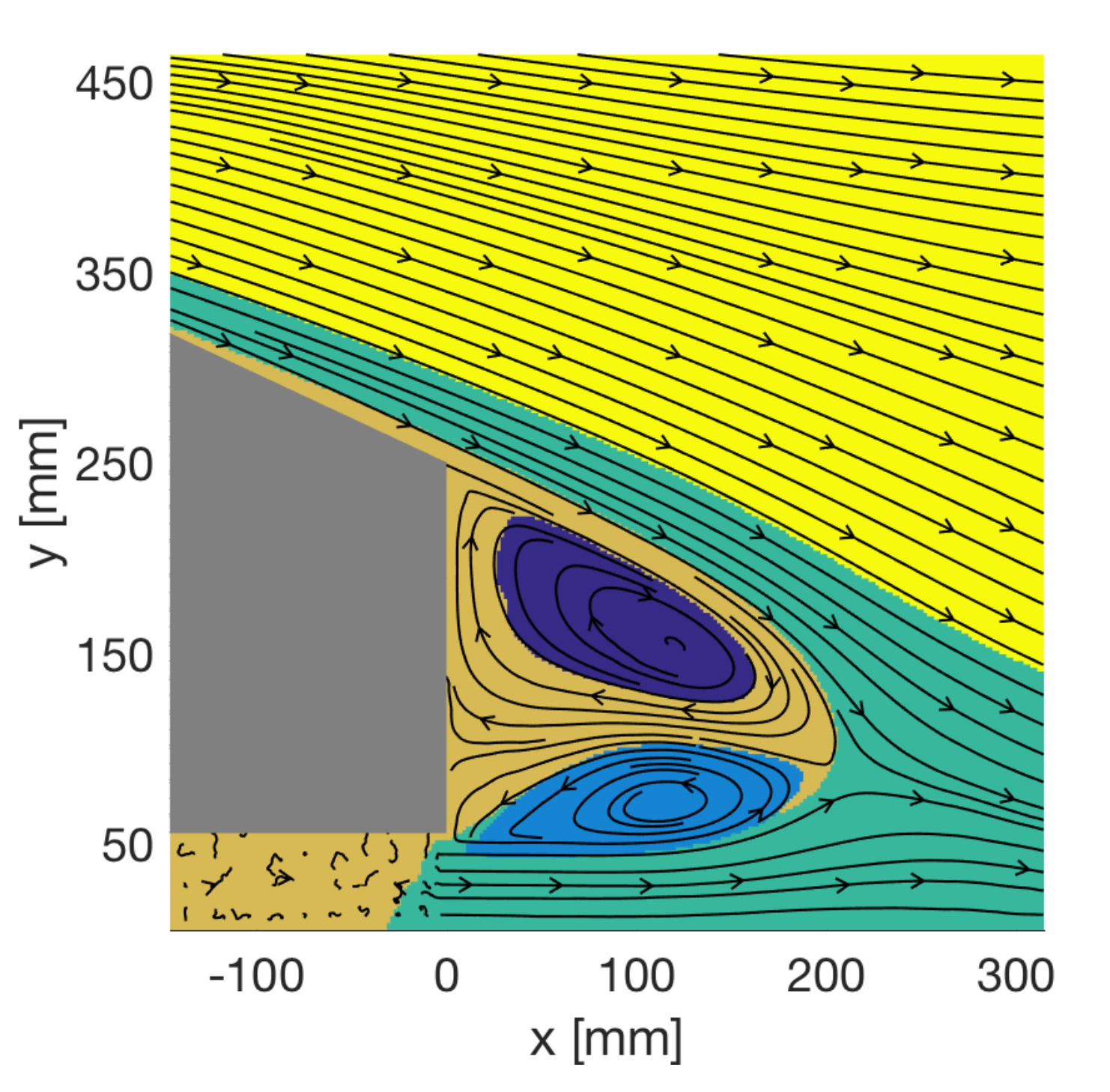}
       \caption{5 clusters produced with k-means clustering w.r.t to the eigenvectors from Figure~\ref{fig:chosenEv_side}.}
        \label{fig:cluster_side}
    \end{subfigure}
 \caption{LKDmap analysis of Ahmed body.}
\end{figure}

We consider a data set obtained in a wind tunnel experiment focused on a flow around the rear end of a simplified vehicle model, the so-called modified Ahmed body, in~\cite{rossitto2016influence}. The data set consists of $65\,025$ points in the $xy$ plane with $x$ being the wind direction and $y$ the direction orthogonal to the surface. The corresponding velocities are obtained from PIV\footnote{Particle image velocimetry (PIV) is an optical method which produces two or three-dimensional vector fields of the flow. Temporal dependence of the velocity field is averaged out.} measurements. We construct the LKDmap matrix $L_{\ep, \tilde \ep}$ with $b(x_i) = v(x_i)$ and $A = \beta^{-1} I$. We choose the parameter values $\beta^{-1} = 0.02\,\textup{m}^{2}\textup{s}^{-1}$, $\ep = 0.2\, \textup{ms}$, and $\tilde \ep =0.01\,\textup{mm}^2$.  Figure~\ref{fig:chosenEv_side} shows the eigenvectors $\psi_1$, $\psi_4$, $\psi_6$, and $\psi_{10}$ of $L_{\ep,\tilde\ep}$ together with a visualization of the flow lines. The rear of the Ahmed body is shown in grey. Two counterrotating vortices form directly behind it. There is a clear separation between the recirculation zone containing the vortices and free flow at the top and bottom. The location of this boundary, which is of interest in applications, is clearly identified by $\psi_4$. The other shown eigenfunctions highlight other important aspects of the flow field. The eigenfunctions not shown further stratify the free stream region into thinner layers which are hard to cross in the finite time it takes for particles to travel the open measurement domain. We also note that velocity measurements at the Ahmed body are very noisy and not accessible by the measurement technique below the model, hence randomly assigned in the data postprocessing\footnote{Usually, this part is cut out because only the flow around the car model is relevant. In our analysis, we do not remove this part of the data.}.

Figure~\ref{fig:cluster_side} shows a clustering of the domain into 5 clusters using the popular k-means method on the eigenfunctions $\psi_1, \psi_4, \psi_6$, and $\psi_{10}$, which is equivalent to spectral clustering \cite{shi2000normalized}. In Figure~\ref{fig:cluster_side}, the clustering allows to find the most important structures in the flow: the two main vortices (dark and azure blue) with the surrounding recirculation zone (orange), the free flow (yellow) and the mixing layer (green).

We note that the velocity field is treated here as being stationary. Time dependent velocity fields require an analysis based on coherent rather then almost invariant structures \cite{BaKo16}, but they also can be based on generators of certain space-time processes~\cite{FrKo17}.

\section{Conclusion}

In this article, we have extended the original diffusion maps construction to a much larger class of differential operators that can be approximated. The differential operators that we approximate are the generators of certain dynamical It\^o diffusion processes that live on the data manifold $\op{M}$. Our approximations are geometric; it is not necessary that the data was generated by the dynamical process considered, it only has to sample $\op{M}$. We stressed the importance of differentiating between geometry and sampling, and all our approximations are asymptotically unbiased in the sense that they do not depend on the sampling density $q$ in the limit where the scale parameter $\ep$ tends to zero.

Our first extension, TMDmap, allows users to apply the popular diffusion map construction with the $\alpha = 1/2$ normalization \cite{Nadler2006, Nadler2008, singer2009detecting, szlam2008regularization} in cases where sampling from the density of interest is too difficult. This is the case in many applications where the density of interest is multimodal, and TMDmap allows the full toolbox of importance sampling and enhanced sampling strategies to be incorporated into diffusion maps. Sampling strategies often have a bias-variance trade off, and with TMDmap one may leverage a low variance and high bias sampler by removing the bias at the stage of the diffusion map construction. We leave further exploration of this idea for future work.

Our second extension, LKDmap, generalizes TMDmap further and allows the approximation of the forward and backward Fokker--Planck operators corresponding to a large class of It\^o diffusions, including non-gradient flow and anisotropic noise. We draw from the local kernels theory developed in \cite{BeSa16}, and the class of It\^o diffusions that can be approximated with LKDmap is equal to the class corresponding to the prototypical local kernels in \cite{BeSa16}. We cannot handle degenerate noise, which would render the diffusion matrix non positive definite. Two application areas come to mind for LKDmap. First, if the data points are given in conjunction with velocity vectors, then LKDmap can approximate the generator of the flow corresponding to the velocity field with added diffusive noise. We showed in Section \ref{sec:numerics_d} that this leads to the identification of almost invariant structures in the flow. Similarly, having sampled the significant part of the state space of a process, one may wish to analyze how the global dynamical properties---connected to dominant eigenmodes of the generator---change under slight variations of the dynamics. Molecular dynamics comes to mind as an immediate application field; biasing the molecular potential to influence the statistical weights of certain conformations and varying dominant implied time scales.
Second, one may wish to approximate high-dimensional dynamics by an It\^o diffusion in a reduced space. In Section \ref{sec:numerics_c}, we showed that LKDmap together with estimators for the drift and diffusion coefficients can approximate the generator of such an effective dynamics and is able to uncover a hidden slow variable that is a nonlinear function of the data.

Both TMDmap and LKDmap have scale parameters that need to be selected, which can be difficult in practice. The single scale parameter $\ep$ in TMDmap plays the same role as $\ep$ in the original diffusion maps construction \cite{CoLa06}, and similar considerations for selecting it apply. The optimal choice of $\ep$ in diffusion maps has been discussed in \cite{SINGER2006}, but it involves prefactors that depend on the manifold $\op{M}$ and are unknown in practice. The connectivity of the neighbourhood graph and the sparsity of the kernel matrix might offer good practical indicators for choosing $\ep$. A multiscale analysis such as in \cite{Little2009} may give more insight at a higher computational cost. LKDmap has two scale parameters $\ep$ and $\tilde \ep$ with different interpretations. The parameter $\tilde\ep$, which is used to construct the kernel density estimate of $q$ and should be thought of as a length scale, may be selected purely based on geometric considerations. The second scale $\ep$, which should rather be thought of as a time scale, may be selected based on typical time scales of the dynamics which is to be approximated.

A central aspect of this work is the usage of a kernel density estimator for unbiasing. We only use the most simple kernel density estimators based on radially symmetric kernels with fixed bandwidth, but in principle one could use more sophisticated ones based on e.g. variable bandwidth. This may reduce the sensitivity of the algorithms with respect to $\ep$ and it may lead to improved convergence properties. Diffusion maps with variable bandwidth kernels have been considered in \cite{Berry2016}. We leave further considerations for future work.

Another direction for future work is time dependent data, either in the form of time ordered position measurements of a realization of the dynamics, or in the form of time dependent velocity fields. Diffusion maps for changing data has been considered in \cite{Hirn2014} and \cite{BaKo16}.

\section*{Acknowledgements}

Zofia Trstanova gratefully acknowledges funding from the EPSRC Grant EP/P006175/1. The other authors acknowledge funding from the DfG through the collaborative research center CRC 1114 and from Freie Universit\"at Berlin. We would like to thank to Giacomo Rossitto for providing the flow data set.

\begin{appendices}

\numberwithin{equation}{section}

\section{Preliminaries}\label{sec:preliminaries}

Following the lines of \cite{BeSa16}, we introduce zeroth, first and second order moments of the local kernel $k^{A,b}_\ep$:
\begin{align}
m(x) & \; = \lim_{\ep\rightarrow 0} \int_{T_x\op{M}} k_\ep^{A,b}\left(x, x + \sqrt{\ep} \hat z\right) dz\label{eq:kernel_moment0}\\
\mu_i(x) & \; = \lim_{\ep\rightarrow 0} \frac{1}{\sqrt{\ep}}\int_{T_x\op{M}} z_i k_\ep^{A,b}\left(x, x + \sqrt{\ep} \hat z\right) dz\\
C_{ij}(x) & \; = \lim_{\ep\rightarrow 0} \int_{T_x\op{M}}z_i z_j k_\ep^{A,b}\left(x, x + \sqrt{\ep} \hat z\right) dz
\end{align}
where $\hat z$ is $z$ projected onto the tangent space $T_x\op{M}$ of $\op{M}$ at $x\in \op{M}$. It was shown in \cite{BeSa16} that with $k_\ep^{A,b}$ defined via \eqref{eq:localkernel}, we have $\mu(x) = m(x)b(x)$ and $C(x) = m(x) A(x)$. The following result was shown in \cite{BeSa16} by combining Lemmas 3.9 and 3.11 therein. The additional factor of $m(x)$ is due to the fact that the operators $\op{L}$ and $\op{L}^*$ are defined via $\mu$ and $C$ in \cite{BeSa16}; we prefer to define $\op{L}$ and $\op{L}^*$ with $b = m^{-1}\mu$ and $A=m^{-1}C$ instead.
\begin{lemma}\label{lemma:expansions}
Let $A(x)$ be a matrix-valued function on $\op{M}\subset \R^N$ such that each $A(x)$ is a symmetric positive definite $N\times N$ matrix, and let $b(x)$ be a vector-valued function. Let $k_\ep^{A,b}$ be defined via \eqref{eq:localkernel} and let $f : \op{M} \rightarrow \R$ be smooth. We define the integral operators
\begin{subequations}
\begin{align}
\op{G}_\ep f(x) &\; = \ep^{-d/2}\int_\op{M} k^{A,b}_\ep (x,y)f(y) dy \label{eq:G}\\
\op{G}^*_\ep f(x) &\; = \ep^{-d/2}\int_\op{M} k^{A,b}_\ep (y,x)f(y) dy\label{eq:Gstar}.
\end{align}
\end{subequations}
Then the expansions
\begin{subequations}
\begin{align}
\op{G}_\ep f(x) &\; = m(x)f(x) + \ep\left[\omega(x)f(x) + m(x)\op{L}f(x)\right] + \mathcal{O}(\ep^{3/2}), \label{eq:G_expand}\\
\op{G}^*_\ep f(x) &\; = m(x)f(x) + \ep\left[(\omega(x)f(x) + \op{L}^*(mf)(x)\right] + \mathcal{O}(\ep^{3/2})\label{eq:Gstar_expand}
\end{align}
\end{subequations}
hold, where $\omega(x)$ depends on the kernel and induced metric $g$, and $m(x)$ is the zeroth moment of $k_\ep^{A,b}$ defined in \eqref{eq:kernel_moment0}.
\end{lemma}

Since the data points $x_i\in \set{D}^{(m)}$ are sampled according to the density $q$, the action of matrices on vectors defined on the data converges to integrals with respect to $q$ in the Monte Carlo sense. For example,
\[
\lim_{m\rightarrow\infty} \frac{1}{m} (K^{A,b}_\ep [f])_k = \int_\op{M} k^{A,b}_\ep(x_k,y) q(y) f(y) dy = \ep^{d/2} \op{G}_\ep (qf)(x_k).
\]
For a finite value of $m$, the relative error in this expression is expected to be of order $\mathcal{O}(m^{-1/2}\ep^{-d/4})$ \cite{CoLa06}. It turns out that one can use cancellations in order to improve the $\ep$-order in the final variance terms for the matrices of interest \cite{SINGER2006}.

Let $k_\ep$ again be the isotropic kernel defined in \eqref{eq:kernel}. For the purposes of this section, it will be convenient to distinguish between the \emph{kernel density estimate}, which we rename here to
\[
\hat q_\ep(x_i) = \sum_{j=1}^m k_\ep(x_i, x_j),
\]
and the (properly normalized) convolution of $q$ with the kernel $k_\ep$
\begin{equation}
\label{eq:q_convolution}
q_\ep(x) = \ep^{-d/2}\int_\op{M} k_\ep(x,y) q(y) dy.
\end{equation}
Then $\lim_{m\rightarrow \infty} m^{-1} \hat q_\ep(x_i) = \ep^{d/2} q_\ep(x_i)$. The following result was shown in \cite{CoLa06}. It can also be deduced from Lemma \ref{lemma:expansions} by noting that the kernel $k_\ep$ is equal to $k_\ep^{A,b}$ with $A = I$ and $b = 0$, and the zeroth moment of $k_\ep$ is $m(x) = \ep^{d/2}$.

\begin{lemma}
\label{lemma:q}
The function $q_\ep$ defined in \eqref{eq:q_convolution} has an expansion of the form
\[
q_\ep = q + \ep(\omega q - \Delta q) + \mathcal{O}(\ep^{3/2})
\]
uniformly on $\op{M}$. Sufficiently far away from $\partial \op{M}$, the order of the error term is improved to $\mathcal{O}(\ep^2)$.
\end{lemma}

\section{Proofs}\label{sec:proofs}

\subsection{Proof of Theorem \ref{thm:target_diffmap}}\label{sec:proof_TMDmap}

We first consider the matrix $K_{\ep,\pi} = K_\ep D_{\ep,\pi}$ with $(D_{\ep,\pi})_{ii} = \pi^{1/2}(x_i)\hat q^{-1}_\ep(x_i)$. Note that
\[
\lim_{m\rightarrow\infty} \frac{1}{m} (K_{\ep,\pi} [f])_k  = \lim_{m\rightarrow\infty}\frac{1}{m}\sum_{j=1}^m \frac{k_\ep(x_k,x_j)}{\hat q_\ep(x_j)}\pi^{1/2}(x_j)f(x_j).
\]
With the notation in \eqref{eq:q_convolution} and Lemma \ref{lemma:expansions}, we thus have, with $\tilde q_\ep = \pi^{1/2} q q_\ep^{-1}$,
\[
\lim_{m\rightarrow\infty} \frac{1}{m} (K_{\ep,\pi} [f])_k = \ep^{-d/2}\int k_\ep(x_k,y) f(y) \frac{\pi^{1/2}(y)q(y)}{q_\ep(y)} dy = \op{G}_\ep (\tilde q_\ep f)(x_k).
\]
The operator $\op{L}$ In the corresponding expansion \eqref{eq:G_expand} of $\op{G}_\ep$ is given by $\op{L} = \Delta$. The limiting expression for $L_{\ep,\pi}$ as $m\rightarrow \infty$ reads
\begin{equation}
\label{eq:proof_TMDmap_1}
\lim_{m\rightarrow \infty} (L_{\ep,\pi}[f])_k = \lim_{m\rightarrow \infty} \ep^{-1} \left( \frac{(K_{\ep,\pi}[f])_k}{(K_{\ep,\pi}[\mathbf{1}])_k} - f(x_k)\right) = \ep^{-1} \left( \frac{\op{G}_\ep( \tilde q_\ep f)(x_k)}{\op{G}_\ep(\tilde q_\ep)(x_k)} - f(x_k)\right).
\end{equation}
Lemma \ref{lemma:q} with the shorthand $q^{(1)} = \omega - q^{-1}\Delta q$ implies $\tilde q_\ep = \pi^{1/2}(1 - \ep q^{(1)}) + \mathcal{O}(\ep^{3/2})$ uniformly on $\op{M}$. We invoke \eqref{eq:G_expand} from Lemma \ref{lemma:expansions} to obtain
\begin{align*}
\op{G}_\ep(\tilde q_\ep f) & \; = m f \tilde q_\ep + \ep\left[\omega\tilde q_\ep f + m \Delta(\tilde q_\ep f)\right] + \mathcal{O}(\ep^{3/2}) \\
&\; = m\pi^{1/2}f + \ep \pi^{1/2}\left[\omega f - q^{(1)}m f + m\pi^{-1/2}\Delta(\pi^{1/2} f)\right] + \mathcal{O}(\ep^{3/2})
\end{align*}
where $m(x)$ is the zeroth moment of the kernel $k_\ep$, defined in \eqref{eq:kernel_moment0}. Consequently,
\[
\op{G}_\ep \tilde q_\ep = m\pi^{1/2} + \ep \pi^{1/2}\left[\omega  - q^{(1)}m  + m\pi^{-1/2}\Delta(\pi^{1/2}) \right] + \mathcal{O}(\ep^{3/2}).
\]
Dividing the two equations gives, with the shorthand $\tilde \pi = \pi^{1/2}$,
\begin{align}
\left(\op{G}_\ep \tilde q_\ep \right)^{-1}\op{G}_\ep(\tilde q_\ep f) =\;& f + \ep\left[ m^{-1}\omega f- q^{(1)}f + \tilde\pi^{-1}\Delta(\tilde\pi f)\right] \notag \\
&\; - \ep\left[ m^{-1}\omega f - q^{(1)}f + f\tilde\pi^{-1}\Delta\tilde\pi\right] + \mathcal{O}(\ep^{3/2}) \notag\\
=\;& f + \ep\left[\tilde\pi^{-1}\Delta(\tilde\pi f) - f\tilde\pi^{-1}\Delta\tilde\pi\right]+ \mathcal{O}(\ep^{3/2}) \notag\\
=\; & f + \ep\op{L} f+ \mathcal{O}(\ep^{3/2}), \label{eq:proof_TMDmap_2}
\end{align}
where the last line follows from
\[
\tilde\pi^{-1}\Delta(\tilde\pi f) - f\tilde\pi^{-1}\Delta\tilde\pi = \Delta f + 2 \tilde\pi^{-1}\nabla\tilde\pi \cdot\nabla f= \Delta f + 2 \nabla(\log\tilde\pi) \cdot\nabla f = \Delta f +  \nabla(\log\pi) \cdot\nabla f = \op{L}f. 
\]
Theorem \ref{thm:target_diffmap} now follows from \eqref{eq:proof_TMDmap_1} and \eqref{eq:proof_TMDmap_2}. Sufficiently far away from $\partial \op{M}$ the order of the error term improves to $\mathcal{O}(\ep^2)$. $\qed$

\subsection{Proof of Theorem \ref{thm:localkernels_q}}\label{sec:proof_localkernel}

With the notation in \eqref{eq:q_convolution} and Lemma \ref{lemma:expansions}, we have
\begin{equation}
\label{eq:proof_localkernels_1}
\lim_{m\rightarrow \infty}\frac{1}{m} K_\ep^{A,b} [f]_k = \int_{\op{M}} k^{A,b}_\ep(x_k, y) q(y) f(y) dy = \ep^{d/2}\op{G}_\ep(qf)(x_k).
\end{equation}
The limiting expression for $L_{\ep}$, defined in \eqref{eq:localkernels_L}, as $m\rightarrow \infty$ then reads
\begin{equation*}
\lim_{m\rightarrow \infty} (L_{\ep}[f])_k = \lim_{m\rightarrow \infty} \ep^{-1} \left( \frac{(K_{\ep}^{A,b} [f])_k}{(K_{\ep}^{A,b} [\mathbf{1}])_k} - f(x_k)\right) = \ep^{-1} \left( \frac{\op{G}_\ep( q f)(x_k)}{\op{G}_\ep(q)(x_k)} - f(x_k)\right).
\end{equation*}
We invoke \eqref{eq:G_expand} from Lemma \ref{lemma:expansions} to obtain
\[
\op{G}_\ep(q f)  = m f q + \ep\left[\omega q f + m \op{L}(q f)\right] + \mathcal{O}(\ep^{3/2})
\]
where, with abuse of notation, $m = m(x)$ denotes the zeroth moment of the kernel $k_\ep^{A,b}$, defined in \eqref{eq:kernel_moment0}. Consequently
\begin{equation}
\label{eq:proof_localkernels_1a}
\left(\op{G}_\ep q\right)^{-1}\op{G}_\ep(q f) = f + \ep \left[q^{-1}\op{L}(qf) - q^{-1}f\op{L}q\right]   +  \mathcal{O}(\ep^{3/2}).
\end{equation}
The first equation in Theorem \ref{thm:localkernels_q} now follows from \eqref{eq:backward_gen} and
\begin{align*}
q^{-1}\op{L}(qf) - q^{-1}f\op{L}q & \; = q^{-1} \left[b\cdot \nabla (qf) + A_{ij}\nabla_i\nabla_j (qf)\right] - q^{-1}f[b\cdot \nabla q + \nabla_i\nabla_j q] \\
& \; = \op{L}f + 2q^{-1}A_{ij}(\nabla_i f)(\nabla_j q).
\end{align*}
For the dual result, we note that, from \eqref{eq:proof_localkernels_1} and \eqref{eq:localkernels_Lstar},
\begin{equation}
\label{eq:proof_localkernels_2}
\lim_{m\rightarrow \infty} (L^*_{\ep}[f])_k = \ep^{-1} \left(\op{G}^*_\ep\left((\op{G}_\ep q)^{-1}q f\right)(x_k) - f(x_k)\right).
\end{equation}
We invoke \eqref{eq:Gstar_expand} from Lemma \ref{lemma:expansions} to obtain
\[
\op{G}^*_\ep(q f)  = m f q + \ep\left[\omega q f + \op{L}^*(m q f)\right] + \mathcal{O}(\ep^{3/2}).
\]
Consequently
\begin{equation*}
\op{G}^*_\ep\left((\op{G}_\ep q)^{-1}q f\right) = f + \ep \left[\op{L}^*f - q^{-1}f\op{L}q\right]   +  \mathcal{O}(\ep^{3/2}).
\end{equation*}
The second equation in Theorem \ref{thm:localkernels_q} now follows from \eqref{eq:forward_gen}. $\qed$

\subsection{Proof of Theorem \ref{thm:local_diffmap}}\label{sec:proof_LKDmap}

We first consider the matrix $\tilde K_{\ep,\tilde \ep}^{A,b} = K_\ep^{A,b} D_{\tilde\ep}$ with $(D_{\tilde\ep})_{ii} = \hat q_{\tilde\ep}(x_i)^{-1}$. Note that
\[
\lim_{m\rightarrow\infty} \frac{1}{m} (\tilde K_{\ep,\tilde \ep}^{A,b} [f])_k  = \lim_{m\rightarrow\infty}\frac{1}{m}\sum_{j=1}^m \frac{k^{A,b}_\ep(x_k,x_j)}{\hat q_{\tilde\ep}(x_j)}f(x_j).
\]
With the notation in \eqref{eq:q_convolution} and Lemma \ref{lemma:expansions}, we thus have, with $\tilde q_{\tilde \ep} = q q_{\tilde \ep}^{-1}$,
\[
\lim_{m\rightarrow\infty} \frac{1}{m} (\tilde K_{\ep,\tilde \ep}^{A,b} [f])_k = (\tilde\ep / \ep)^{-d/2}\int k^{A,b}_\ep(x_k,y)\frac{q(y)}{q_{\tilde \ep}(y)}f(y)  dy =  (\tilde\ep / \ep)^{-d/2}\op{G}_\ep (\tilde q_{\tilde \ep}f) (x_k).
\]
The limiting expression for $L_{\ep,\tilde\ep}$ as $m\rightarrow \infty$ reads
\begin{equation}
\label{eq:proof_LKDmap_1}
\lim_{m\rightarrow \infty} (L_{\ep,\tilde \ep}[f])_k = \lim_{m\rightarrow \infty} \ep^{-1} \left( \frac{(\tilde K_{\ep,\tilde \ep}^{A,b} [f])_k}{(\tilde K_{\ep,\tilde \ep}^{A,b} [\mathbf{1}])_k} - f(x_k)\right) = \ep^{-1} \left( \frac{\op{G}_\ep( \tilde q_\ep f)(x_k)}{\op{G}_\ep(\tilde q_\ep)(x_k)} - f(x_k)\right).
\end{equation}
Lemma \ref{lemma:q} with the shorthand $q^{(1)} = \omega - q^{-1}\Delta q$ implies $\tilde q_{\tilde\ep} = 1 - \tilde\ep q^{(1)} + \mathcal{O}(\tilde\ep^{3/2})$ uniformly on $\op{M}$. We invoke \eqref{eq:G_expand} from Lemma \ref{lemma:expansions} to obtain
\begin{align*}
\op{G}_\ep(\tilde q_{\tilde\ep} f) & \; = m f \tilde q_{\tilde\ep} + \ep\left[\omega\tilde q_{\tilde\ep} f + m \op{L}(\tilde q_{\tilde\ep} f)\right] + \mathcal{O}(\ep^{3/2}) \\
&\; = m f + \ep \left[\omega f  + m\op{L}f\right] - \tilde \ep q^{(1)}m f + \mathcal{O}(\ep, \tilde \ep) + \mathcal{O}(\ep^{3/2}) + \mathcal{O}(\tilde \ep^{3/2})
\end{align*}
where, with abuse of notation, $m = m(x)$ denotes the zeroth moment of the kernel $k_\ep^{A,b}$, defined in \eqref{eq:kernel_moment0}. Consequently,
\begin{equation}
\label{eq:proof_LKDmap_1a}
\op{G}_\ep(\tilde q_{\tilde\ep}) = m + \ep \omega - \tilde \ep q^{(1)}m  + \mathcal{O}(\ep, \tilde \ep) + \mathcal{O}(\ep^{3/2}) + \mathcal{O}(\tilde \ep^{3/2}).
\end{equation}
Dividing the two equations gives
\begin{align}
\left(\op{G}_\ep \tilde q_{\tilde\ep} \right)^{-1}\op{G}_\ep(\tilde q_{\tilde\ep} f) =\;& f + \ep\left[ m^{-1}\omega f + \op{L}f\right] - \tilde \ep q^{(1)}f  \notag \\
&\; - \ep m^{-1}\omega f +\tilde\ep q^{(1)}f + \mathcal{O}(\ep, \tilde \ep) + \mathcal{O}(\ep^{3/2}) + \mathcal{O}(\tilde \ep^{3/2}) \notag\\
=\; & f + \ep\op{L} f + \mathcal{O}(\ep, \tilde \ep) + \mathcal{O}(\ep^{3/2}) + \mathcal{O}(\tilde \ep^{3/2}). \label{eq:proof_LKDmap_2}
\end{align}
The first equation in Theorem \ref{thm:local_diffmap} now follows from combining \eqref{eq:proof_LKDmap_1} and \eqref{eq:proof_LKDmap_2}. Sufficiently far away from $\partial \op{M}$ the order of the last two error terms improves to $\mathcal{O}(\ep^2)$ and $\mathcal{O}(\tilde\ep^2)$ respectively.

For the dual result, note that
\[
\lim_{m\rightarrow\infty} (L_{\ep,\tilde \ep}^* [f])_k  = \lim_{m\rightarrow\infty}\sum_{j=1}^m \frac{k^{A,b}_\ep(x_j,x_k)}{\hat q_{\tilde\ep}(x_j)}(\tilde D_{\ep,\tilde \ep})_{jj}^{-1}f(x_j).
\]
On the other hand, by the definition of $\tilde D_{\ep,\tilde\ep}$ and with $\tilde q_{\tilde \ep} = q q_{\tilde \ep}^{-1}$ as above,
\[
\lim_{m\rightarrow\infty} \frac{1}{m}(\tilde D_{\ep,\tilde\ep})_{jj} = \lim_{m\rightarrow\infty} \tilde K_{\ep,\tilde \ep}^{A,b} [\mathbf{1}])_j = \op{G}_\ep \tilde q_{\tilde\ep} (x_j).
\]
Combining the two equations gives the limiting expression for $L^*_{\ep,\tilde\ep}$ as $m\rightarrow \infty$ as
\begin{equation*}
\lim_{m\rightarrow \infty} (L^*_{\ep,\tilde \ep}[f])_k = \ep^{-1} \left(\op{G}^*_\ep\left((\op{G}_\ep\tilde q_{\tilde\ep})^{-1}\tilde q_{\tilde\ep} f\right)(x_k) - f(x_k)\right).
\end{equation*}
Let $g = (\op{G}_\ep\tilde q_{\tilde\ep})^{-1} f$. Using \eqref{eq:Gstar_expand} and $\tilde q_{\tilde\ep} = 1 - \tilde\ep q^{(1)} + \mathcal{O}(\tilde\ep^{3/2})$ gives (again with $m = m(x)$ being the zeroth moment \eqref{eq:kernel_moment0})
\begin{align*}
\label{eq:Gstar_expand2}
\op{G}^*_\ep (\tilde q_{\tilde\ep}g) &\; = m \tilde q_{\tilde\ep}g + \ep\left(\omega \tilde q_{\tilde\ep}g  + \op{L}^*(m \tilde q_{\tilde\ep}g)\right) + \mathcal{O}(\ep^{3/2})\\
&\; = m g + \ep\left(\omega g  + \op{L}^*(m g)\right) - \tilde \ep mq^{(1)}g + \mathcal{O}(\ep, \tilde \ep) + \mathcal{O}(\ep^{3/2}) + \mathcal{O}(\tilde \ep^{3/2}).
\end{align*}
On the other hand, using \eqref{eq:proof_LKDmap_1a} gives
\[
g = (\op{G}_\ep\tilde q_{\tilde\ep})^{-1} f = m^{-1}f - \ep m^{-1} \omega f + \tilde\ep q^{(1)} f + \mathcal{O}(\ep, \tilde \ep) + \mathcal{O}(\ep^{3/2}) + \mathcal{O}(\tilde \ep^{3/2}).
\]
Combining the two last equations:
\begin{align*}
\op{G}_\ep^*\left((\op{G}_\ep\tilde q_{\tilde\ep})^{-1}\tilde q_{\tilde\ep} f\right) =\;& f + \ep\left[ m^{-1}\omega f + \op{L}^* f \right] - \tilde \ep q^{(1)}f \\
&\; - \ep m^{-1}\omega f  - \tilde \ep q^{(1)}f + \mathcal{O}(\ep, \tilde \ep) + \mathcal{O}(\ep^{3/2}) + \mathcal{O}(\tilde \ep^{3/2})\\
=\;& f + \ep \op{L}^*f + \mathcal{O}(\ep, \tilde \ep) + \mathcal{O}(\ep^{3/2}) + \mathcal{O}(\tilde \ep^{3/2}).
\end{align*}
The second equation in Theorem \ref{thm:local_diffmap} now follows. Sufficiently far away from $\partial \op{M}$ the order of the last two error terms improves to $\mathcal{O}(\ep^2)$ and $\mathcal{O}(\tilde\ep^2)$ respectively. $\qed$

\end{appendices}

\bibliographystyle{plain}
\bibliography{References,MyBibliography}

\end{document}